\documentstyle[11pt]{article}
\begin{document}
\author{{Shiqiu Zheng$^{1, 2}$\thanks{E-mail: shiqiumath@163.com (S. Zheng).}\ , \ \ Shoumei Li$^{1}$\thanks{Corresponding author, E-mail: li.shoumei@outlook.com (S. Li).}}
  \\
\small(1, College of Applied Sciences, Beijing University of Technology, Beijing 100124, China)\\
\small(2, College of Sciences, North China University of Science and Technology, Tangshan 063009, China)\\
}
\date{}
\title{\textbf{On the representation for dynamically consistent nonlinear evaluations: uniformly continuous case}}\maketitle

\textbf{Abstract:}\ \  A system of dynamically consistent nonlinear evaluation (${\cal{F}}$-evaluation)
provides an ideal characterization for the dynamical behaviors of risk measures
and the pricing of contingent claims. The purpose of this paper is to study the
representation for the ${\cal{F}}$-evaluation by the solution of a backward stochastic differential equation (BSDE).
Under a general domination condition, we prove that any ${\cal{F}}$-evaluation can be represented by the solution of a BSDE with a generator which is Lipschitz in $y$ and uniformly continuous in $z$. \\



\section{Introduction}
The notion of $g$-expectation was introduced by Peng [12] in 1997 via the
solution of a BSDE. $g$-expectation is a dynamically consistent nonlinear expectation and it has many applications in utilities and risk measures.
An axiomatic system of dynamically consistent nonlinear expectation
(${\cal{F}}$-expectation for short) was introduced by Coquet et al. [3] in 2002.
It was shown in Coquet et al. [3] that, under a certain domination condition,
any ${\cal{F}}$-expectation can be represented as a $g$-expectation.
Note that the $g$-expectation involved in the representation theorem in [3]
was defined by a BSDE with a generator Lipschitz in $z$ and independent of $y$.
As an extension of the representation in [3] to a L\'{e}vy filtration,
Royer [18] obtained a corresponding representation by $g$-expectations
defined via BSDEs with jumps.
In 2012, Cohen [2] obtained a corresponding representation, within a
general filtration, by a $g$-expectation defined via a BSDE in a
general probability space.
Note that the domination conditions in [18] and [2] are both similar to that of [3].
Consequently, the $g$-expectations involved in the representation theorems in
[18] and [2] are both defined by BSDEs with Lipschitz generators.
In 2008, Hu et al. [7] considered a quadratic ${\cal{F}}$-expectation, and
showed that, under three domination conditions, any quadratic
${\cal{F}}$-expectation can be represented as a $g$-expectation defined by
a BSDE with a quadratic growth.
Recently, under a domination condition more general than that of [3],
Zheng and Li [19] obtained a representation theorem by a $g$-expectation
defined by a BSDE whose generator is independent of $y$, uniformly continuous in $z$.

It is well known that the famous Black-Scholes option pricing model is a linear BSDE.
Peng [16] defined a $g$-evaluation by the solution of a nonlinear BSDE
and used this $g$-evaluation as a general pricing model.
For quadratic $g$-evaluations, we refer to Ma and Yao [11].
In [14, 16], Peng introduced an axiomatic system of dynamically consistent nonlinear evaluation (${\cal{F}}$-evaluation for short).
The concept of ${\cal{F}}$-evaluation is a natural extension of the concept of
${\cal{F}}$-expectation.
In [13, 14], Peng proved that any ${\cal{F}}$-evaluation ${\cal{E}}_{s,t}[\cdot]$ is a $g$-evaluation under the following domination condition:
$$
{\cal{E}}_{s,t}[X]-{\cal{E}}_{s,t}[Y]\leq{\cal{E}}_{s,t}^{\mu,\mu}[X-Y],\eqno(1.1)
$$
where ${\cal{E}}_{s,t}^{\mu,\mu}[\cdot]$ is a $g$-evaluation defined by the
solution of a BSDE with a generator of the form $g(y,z)=\mu|y|+\mu|z|$ for
some constant $\mu>0.$
Note that the $g$-evaluation involved in the representation theorem in [14] is defined by a BSDE whose generator is Lipschitz in $y$ and $z$.
Recently, based on the representation in [14], Hu [6] obtained a
representation for ${\cal{F}}$-evaluations with $L^p$ terminal variables
($p>1$) under the domination condition (1.1).

The main reason for studying this topic is that the axiomatic systems of ${\cal{F}}$-evaluations and ${\cal{F}}$-expectations provide an ideal characterization of the dynamical behaviors of
risk measures and the pricing of contingent claims (see [14, 16]).
The representation theorems for ${\cal{F}}$-evaluations and ${\cal{F}}$-expectations
mean that any risk measure and the pricing of contingent claims
can be represented as the solution of  a BSDE under some conditions.
Peng raised the following interesting question in [14]: are the notions of
$g$-expectations and $g$-evaluations general enough to represent all
"regular enough" dynamically consistent nonlinear expectations and evaluations?
This paper is devoted to this question. We will show that
any ${\cal{F}}$-evaluation ${\cal{E}}_{s,t}[\cdot]$ is a $g$-evaluation, provided
that the following general domination condition holds:
$$
{\cal{E}}_{s,t}[X]-{\cal{E}}_{s,t}[Y]\leq{\cal{E}}_{s,t}^{\mu,\phi}[X-Y],\eqno(1.2)
$$
where ${\cal{E}}_{s,t}^{\mu,\phi}[\cdot]$ is a $g$-evaluation defined by the solution of
a BSDE with a generator of the form $g(y,z)=\mu|y|+\phi(|z|),$ where $\mu>0$ is a constant and $\phi(\cdot):{\mathbf{R_+}}\rightarrow{\mathbf{R_+}},$ is a continuous, increasing, subadditive function with $\phi(0)=0$ and satisfies a linear growth condition.
The $g$-evaluation in our representation theorem is defined by a BSDE whose generator is Lipschitz in $y$ and uniformly continuous in $z.$

The main result of this paper is an extension of the main results in [3, 13, 14].
It also generalizes the main result in our recent work [19]. The paper [19]
used a method developed in [3] and heavily dependent on the translation invariance
of the ${\cal{F}}$-expectation.
The present paper follows the methods developed in [14], but the argument is
by no means easy.
For example, some of the estimates crucial in the proof of the main result of  [14] are not true in our setting.
We have to develop some techniques to overcome the various difficulties
arising from the lack of Lipschitz continuity.
Estimate on the solution of the BSDEs and the localization technique
play important roles in our proofs.
We now point out a few differences between the present work and [14].
\begin{itemize}
  \item[(i)]
In [14], the introduction of ${\cal{E}}_{s,t}[\cdot;K]$ needs some
convergence results which guaranteed by the estimates in
[14, Theorem 4.1 and Corollary 5.8].  We establish these convergence results
in our setting by using an approximation method.
We also use a different method to prove the ${\cal{E}}_{s,t}[\cdot]$ admits
a RCLL version (see Lemma 3.11).
  \item[(ii)]
In [14], the definition of ${\cal{E}}_{\sigma,\tau}[\cdot]$ with $\sigma, \tau\in{\cal{T}}_{0,T}$ and the proof of optional stopping theorem for ${\cal{E}}_{s,t}[\cdot]$-supermartingales depend on some $L^2$ estimates given in
[14, Corollary 10.15 and Lemma 10.16].
In this paper, we prove a crucial estimate for ${\cal{E}}^g_{s,t}[\cdot;K]$ in the $L^\infty$ sense for bounded terminal variables and bounded $K$ of the form $K_t=\int_0^t\gamma_sds$ (see Lemma 2.6).
This estimate helps us to get some important convergence results (see Lemma 4.2).
With the help of these convergence results, we extend the definition of ${\cal{F}}$-evaluation ${\cal{E}}_{s,t}[\cdot;K]$  to ${\cal{E}}_{\sigma,\tau}[\cdot;K]$ with $\sigma, \tau\in{\cal{T}}_{0,T}$ for a special kind of $K.$
Moreover, we also prove an optional stopping theorem for locally bounded ${\cal{E}}_{s,t}[\cdot;K]$-supermartingales  (see Lemma 4.7).
  \item[(iii)]
In [14], the fixed point method used to solve the BSDEs under
${\cal{E}}_{s,t}[\cdot]$ depends on the $L^2$ estimate given in [14, Proposition 4.5],  and the Doob-Meyer decomposition is obtained for square integrable ${\cal{E}}_{s,t}[\cdot]$-supermartingales. We solve the BSDEs under ${\cal{E}}_{s,t}[\cdot]$ with bounded terminal variables by using our $L^\infty$ estimate (see Lemma 2.6).  By this and our optional stopping theorem, we prove
 a Doob-Meyer decomposition for locally bounded ${\cal{E}}_{s,t}[\cdot]$-supermartingales  (see Theorem 5.4).

  \item[(iv)]
The proofs of the representation theorems in [2, 3, 14, 18] use a Doob-Meyer decomposition for square integrable ${\cal{E}}_{s,t}[\cdot]$-supermartingales.
The proofs of [7, 19] use a Doob-Meyer decomposition for ${\cal{E}}_{s,t}[\cdot]$-supermartingales with a special structure.
In this paper, we use a localization method based on stopping times to
guarantee that our Doob-Meyer decomposition for locally bounded ${\cal{E}}_{s,t}[\cdot]$-supermartingales is good enough  in our proof.
\end{itemize}

This paper is organized as follows. In the next section, we will recall the definitions of $g$-evaluation, $g$-martingale and prove some important convergence results and estimates. In Section 3, we will recall the definitions of ${\cal{F}}$-evaluation ${\cal{E}}_{s,t}[\cdot]$, ${\cal{E}}_{s,t}[\cdot]$-martingale and prove some useful properties.  In Section 4, we will establish an optional stopping theorem for locally bounded ${\cal{E}}_{s,t}[\cdot;K]$-supermartingales. In Section 5, we will give a Doob-Meyer decomposition for locally bounded ${\cal{E}}_{s,t}[\cdot]$-supermartingales. In Section 6, we will prove the main result of this paper: a representation theorem for ${\cal{F}}$-evaluations.

%%%%%%%%%%%%%%%%%%%%%%%%%%%%%%%%%%%%%%%%%%%%%%%%%%%%%%%%%%%%%%%%
\section{$g$-evaluations and related properties}
%%%%%%%%%%%%%%%%%%%%%%%%%%%%%%%%%%%%%%%%%%%%%%%%%%%%%%%%%%%%%%%%
In this paper, we consider a complete probability space  $(\Omega,\cal{F},\mathit{P})$ on which a $d$-dimensional standard Brownian motion ${{(B_t)}_{t\geq
0}}$ is defined. Let $({\cal{F}}_t)_{t\geq 0}$ denote
the natural filtration generated by ${{(B_t)}_{t\geq 0}}$, augmented
by the $\mathit{P}$-null sets of ${\cal{F}}$. Let $|z|$ denote the
Euclidean norm of $\mathit{z}\in {\mathbf{R}}^d$ and $T>0$ be a given time horizon. For stopping times $\tau_1$ and $\tau_2$ satisfying $\tau_1\leq \tau_2\leq T,$ let ${\cal{T}}_{\tau_1,\tau_2}$ be the set of all stopping times $\tau$ satisfying $\tau_1\leq \tau\leq \tau_2.$ Let ${\cal{T}}^0_{\tau_1,\tau_2}$ be a subset of ${\cal{T}}_{\tau_1,\tau_2}$ such that any member in ${\cal{T}}^0_{\tau_1,\tau_2}$ takes values in a finite set. For $\tau\in{\cal{T}}_{0,T},$ we define the following usual spaces:

$L^2({\mathcal {F}}_\tau;{\mathbf{R}}^d)=\{\xi:{\cal{F}}_\tau$-measurable
${\mathbf{R}}^d$-valued random variable; ${\mathbf{E}}\left[|\xi|^2\right]<\infty\};$

$L^\infty({\mathcal {F}}_\tau;{\mathbf{R}}^d)=\{\xi: {\cal{F}}_\tau$-measurable
${\mathbf{R}}^d$-valued random variable; $\|\xi\|_\infty=\textrm{esssup}_{\omega\in\Omega}|\xi|<\infty\};$

$L^2_{\cal{F}}(0,\tau;{\mathbf{R}}^d)=\{\psi: {\mathbf{R}}^d$-valued predictable
process; $E\left[\int_0^\tau|\psi_t|^2dt\right]
<\infty \};$

$L^\infty_{\cal{F}}(0,\tau;{\mathbf{R}}^d)=\{\psi: {\mathbf{R}}^d$-valued predictable
process; $\|\psi\|_{L^\infty_{\cal{F}}(0,\tau)}=\textrm{esssup}_{(\omega,t)\in\Omega\times [0,\tau]}|\psi _t| <\infty\};$

${\mathcal{D}}^2_{\cal{F}}(0,\tau;{\mathbf{R}}^d)=\{\psi:$ RCLL
process in $L^2_{\cal{F}}(0,\tau;{\mathbf{R}}^d)$;\ $E
[{\mathrm{sup}}_{0\leq t\leq \tau} |\psi _t|^2] <\infty \}$

${\mathcal{D}}^\infty_{\cal{F}}(0,\tau;{\mathbf{R}}^d)=\{\psi: $ RCLL
process in $L^\infty_{\cal{F}}(0,\tau;{\mathbf{R}}^d) \};$

${\mathcal{S}}^2_{\cal{F}}(0,\tau;{\mathbf{R}}^d)=\{\psi:$ continuous
process in ${\mathcal{D}}^2_{\cal{F}}(0,\tau;{\mathbf{R}}^d)\};$

${\mathcal{S}}^\infty_{\cal{F}}(0,\tau;{\mathbf{R}}^d)=\{\psi:$ continuous
process in ${\mathcal{D}}^\infty_{\cal{F}}(0,\tau;{\mathbf{R}}^d)\}.$\\
Note that when $d=1,$ we always denote $L^2({\mathcal {F}}_\tau;{\mathbf{R}}^d)$ by $L^2({\mathcal {F}}_\tau)$ for convention and make the same treatment for the above notations of other spaces.

In this paper, we consider a function $g$
$${g}\left( \omega ,t,y,z\right) : \Omega \times [0,T]\times \mathbf{%
R\times R}^{\mathit{d}}\longmapsto \mathbf{R},$$ such that
$\left(g(t,y,z)\right)_{t\in [0,T]}$ is progressively measurable for
each $(y,z)\in\mathbf{%
R\times R}^{\mathit{d}}$. For the function $g$, in this paper, we make the following assumptions:
\begin{itemize}
  \item (A1). There exists a constant
  $\mu>0$ and a continuous function $\phi(\cdot)$, such that $dP\times dt-a.e.,\
  \forall (y_i,z_i)\in {\mathbf{ R\times R}}^{\mathit{d}},\ \
  (i=1,2):$
  $$|{g}( t,y_{1},z_{1})-{g}( t,y_{2},z_{2}) |\leq \mu|y_{\mathrm 1}
    -y_{2}|+\phi(|z_{\mathrm1}-z_{2}|),$$
  where $\phi(\cdot):{\mathbf{R_+}}\rightarrow{\mathbf{R_+}},$ is subadditive and increasing with $\phi(0)=0$ and has a linear growth with constant $\nu$, i.e., $\forall x\in {\mathbf{R}}^d, \ \phi(|x|)\leq \nu(|x|+1);$
  \item (A2). $\forall (y,z)\in{\mathbf{ R\times R}}^{\mathit{d}},\ g(t,y,z)\in L^2_{\cal{F}}(0,T);$
  \item (A3). $dP\times dt-a.e.,\ g(t,0,0)=0.$
\end{itemize}

For each $(t,y,z)\in[0,T]\times {\mathbf{R}}\times {\mathbf{R}}^d$ and $m>(\mu\vee\nu)$ for $\mu$ and $\nu$ given in (A1), we define
$$\underline{g}_m(t,y,z):=\inf\{g(t,a,b)+m(|y-a|+|z-b|):(a,b)\in{\textbf{Q}}^{1+d}\},\eqno(2.1)$$
$$\overline{g}_m(t,y,z):=\sup\{g(t,a,b)-m(|y-a|+|z-b|):(a,b)\in{\textbf{Q}}^{1+d}\},\eqno(2.2)$$
where $\textbf{Q}$ is the rational set. Note that if $g$ satisfies (A1) and (A2), then by Lepeltier and San Martin [10, Lemma 1], for each $(t,y,z)\in[0,T]\times {\mathbf{R}}\times {\mathbf{R}}^d,$ $\underline{g}_m(t,y,z)$ (resp. $\overline{g}_m(t,y,z)$) is increasing (resp. decreasing) in $m$ and converges to ${g}(t,y,z),$ as $m\rightarrow\infty.$ We also have for each $t\in[0,T],$ $\underline{g}_m(t,y,z)$ (resp. $\overline{g}_m(t,y,z)$) is Lipschitz in $(y,z)$ with constant $m$ and linear growth in $(y,z)$ with constant $(\mu\vee\nu)$.

For $\tau\in{\cal{T}}_{0,T},$ we consider the following BSDE with parameter $(g,\xi,K,\tau):$
$$Y_{\tau\wedge t}=\xi+K_\tau-K_{\tau\wedge t}+\int_{\tau\wedge t}^\tau g\left(s,Y_s,Z_s\right)
ds-\int_{\tau\wedge t}^\tau Z_sdB_s,\ \ \,\ \ t\in[0,T].$$
If the generator $g$ satisfies (A1) and (A2), $\xi\in L^2({\mathcal {F}}_\tau)$ and $K\in {\mathcal{D}}^2_{\cal{F}}(0,T),$ then the BSDE has a unique solution $(Y_t^{g,\xi,K,\tau},Z_t^{g,\xi,K,\tau})\in{\mathcal{D}}^2_{\cal{F}}(0,\tau)\times L^2_{\cal{F}}(0,\tau;{\mathbf{R}}^d)$ (see Jia [8, Theorem 3.6.1]). Furthermore, if $K\in{\mathcal{S}}^2_{\cal{F}}(0,T),$ then $Y_t\in{\mathcal{S}}^2_{\cal{F}}(0,\tau).$ Note that since $\phi$ given in (A1) is subadditive and increasing, then we have $\mu|y|+\phi(|z|)$ satisfies (A1) and (A2). Thus BSDE with parameter $(\mu|y|+\phi(|z|),\xi,K, \tau)$  (resp. $(-\mu|y|-\phi(|z|),\xi,K,\tau)$) has a unique solution.

Now, we introduce the definition of $g$-evaluation, which is introduced by Peng [14, 16] in Lipschitz case, then by Ma and Yao [11] in quadratic case.
\\\\
\textbf{Definition 2.1}  Let $g$ satisfy (A1) and (A2), $K\in {\mathcal{D}}^2_{\cal{F}}(0,T),$ $\sigma,\ \tau\in{\cal{T}}_{0,T}$ and $\sigma\leq\tau.$ Let $\xi\in L^2({\mathcal {F}}_\tau)$ and $(Y_t,Z_t)$ be the solution of BSDE with parameter $(g,\xi,K,\tau)$. We denote the ${\cal{E}}^g_{\sigma,\tau}[\cdot,K]$-evaluation and ${\cal{E}}^g_{\sigma,\tau}[\cdot]$-evaluation of $\xi$ by
$${\cal{E}}^g_{\sigma,\tau}[\xi;K]:=Y_\sigma^{g,\xi,K,\tau},$$
and $${\cal{E}}^g_{\sigma,\tau}[\xi]:={\cal{E}}^g_{\sigma,\tau}[\xi;0].$$

Note that we denote ${\cal{E}}^g_{\sigma,\tau}$ by ${\cal{E}}^{\mu,\phi}_{\sigma,\tau}$ (resp. denote ${\cal{E}}^{g}_{\sigma,\tau}$ by ${\cal{E}}^{-\mu,-\phi}_{\sigma,\tau}$), if $g=\mu|y|+\phi(|z|)$ (resp. $g=-\mu|y|-\phi(|z|)$) for function $\phi(\cdot)$ and constant $\mu>0$, and denote ${\cal{E}}^{g}_{\sigma,\tau}$ by ${\cal{E}}^{\mu,\mu}_{\sigma,\tau}$ (resp. denote ${\cal{E}}^{g}_{\sigma,\tau}$ by ${\cal{E}}^{-\mu,-\mu}_{\sigma,\tau}$),  if $g=\mu|y|+\mu|z|$  (resp. $g=-\mu|y|-\mu|z|$), for constant $\mu>0$. \\

The following Remark 2.2 contains two simple properties of ${\cal{E}}^g_{\sigma,\tau}[\cdot,K]$-evaluations.\\\\
\textbf{Remark 2.2}
Let $g$ satisfy (A1) and (A2), $\sigma, \tau\in{\cal{T}}_{0,T}$ and $\sigma\leq\tau.$ Let $K, K'\in {\mathcal{D}}^2_{\cal{F}}(0,T),$ and $X, X'\in L^2({\mathcal {F}}_\tau).$ Then
\begin{itemize}
\item[(i)]by Jia [8, Theorem 3.6.1], we have
$${\cal{E}}^g_{\sigma,\tau}[X;K]={\cal{E}}^{g^K}_{\sigma,\tau}[X+K_\tau]-K_\sigma,$$
where $g^K(\cdot,\cdot,\cdot):=g(\cdot,\cdot-K_s,\cdot).$
\item[(ii)] by comparison theorem (see Jia [8, Theorem 3.6.3]), we can get
$${\cal{E}}^{-\mu,-\phi}_{\sigma,\tau}[X-X';K-K']\leq{\cal{E}}_{\sigma,\tau}^g[X;K]-{\cal{E}}_{\sigma,\tau}^g[X';K']
\leq{\cal{E}}_{\sigma,\tau}^{\mu,\phi}[X-X';K-K'],$$
from the similar argument as Peng [14, Corollary 4.4].
\end{itemize}
\textbf{Definition 2.3} Let $g$ satisfy (A1) and (A2),  $K\in {\mathcal{D}}^2_{\cal{F}}(0,T) $. A process $Y_t$ with $Y_t\in L^2({\cal{F}}_t)$ for $t\in [0,T],$ is called an ${\cal{E}}^g_{s,t}[\cdot;K]$-martingale (resp. ${\cal{E}}^g_{s,t}[\cdot;K]$-supermartingale, ${\cal{E}}^g_{s,t}[\cdot;K]$-submartingale), if, for each $0\leq s\leq t\leq T,$ we have
$${\cal{E}}^g_{s,t}[Y_t;K]=Y_s,\ \ \ \textmd{ (resp.} \leq,\ \geq\textmd{)}.$$

In the following, we will prove some convergence results and estimates for solutions of BSDEs under (A1) and (A2), which play an important role in this paper. \\\\
\textbf{Lemma 2.5} \textit{Let $g$ satisfy (A1) and (A2), $\tau\in{\cal{T}}_{0,T}.$ Let $K^n, K\in {\mathcal{D}}^2_{\cal{F}}(0,T)$ and $X, X_n\in L^2({\cal{F}}_\tau), n\geq1.$ If $K^n\rightarrow K$ in $L^2_{\cal{F}}(0,T),$ $K^n_\tau\rightarrow K_\tau$ and $X_n\rightarrow X$ both in $L^2({\cal{F}}_T),$ as $n\rightarrow\infty.$ Then we have
$$\lim_{n\rightarrow\infty}E\left[\sup_{s\in[0,T]}|{\cal{E}}_{\tau\wedge s,\tau}^g[X_n;K^n]
+K^n_{\tau\wedge s}-{\cal{E}}_{{\tau\wedge s},\tau}^g[X;K]-K_{\tau\wedge s}|^2\right]=0.$$}\\
\
\emph{Proof.} For $m>(\mu\vee\nu),$ let $\underline{g}_m$ and $\overline{g}_m$ be defined as in (2.1) and (2.2), respectively. Then by comparison theorem (see Jia [8, Theorem 3.6.3]), we have for each $s\in[0,T],$ $${\cal{E}}_{\tau\wedge s,\tau}^{\underline{g}_m}[X_n;K^n]\leq{\cal{E}}_{\tau\wedge s,\tau}^g[X_n;K^n]
\leq{\cal{E}}_{\tau\wedge s,\tau}^{\overline{g}_m}[X_n;K^n],\ \ P-a.s.\eqno(2.3)$$
By Peng [14, Theorem 4.1], we have
$$\lim_{n\rightarrow\infty}E\left[\sup_{s\in[0,T]}|{\cal{E}}_{{\tau\wedge s},\tau}^{\overline{g}_m}[X_n;K^n]
+K^n_{\tau\wedge s}-{\cal{E}}_{ {\tau\wedge s},\tau}^{\overline{g}_m}[X;K]-K_{\tau\wedge s}|^2\right]=0,\eqno(2.4)$$
and
$$\lim_{n\rightarrow\infty}E\left[\sup_{s\in[0,T]}|{\cal{E}}_{{\tau\wedge s},\tau}^{\underline{g}_m}[X_n;K^n]
+K^n_{\tau\wedge s}-{\cal{E}}_{{\tau\wedge s},\tau}^{\underline{g}_m}[X;K]-K_{\tau\wedge s}|^2\right]=0.\eqno(2.5)$$
By (i) in Remark 2.2, the proof of Fan and Jiang [5, Theorem 1] and the uniqueness of solutions, we can get
\begin{eqnarray*}
 \ \ \ \ \ \ \ &&\lim_{m\rightarrow\infty}E\left[\sup_{s\in[0,T]}|{\cal{E}}_{{\tau\wedge s},\tau}^{\overline{g}_m}[X;K]
-{\cal{E}}_{{\tau\wedge s},\tau}^g[X;K]|^2\right]\\
&=&\lim_{m\rightarrow\infty}E\left[\sup_{s\in[0,T]}
|{\cal{E}}_{{\tau\wedge s},\tau}^{\overline{g}_m^K}[X+K_\tau]-{\cal{E}}_{{\tau\wedge s},\tau}^{{g}^K}[X+K_\tau]|^2\right]=0,\ \ \ \ \ \ \ \ \ \ \ \ \ \ \ \ \ \ \ \ \ \ \ \ \ \ \ \  (2.6)
\end{eqnarray*}
and
\begin{eqnarray*}
 \ \ \ \ \ \ \ &&\lim_{m\rightarrow\infty}E\left[\sup_{s\in[0,T]}|{\cal{E}}_{{\tau\wedge s},\tau}^{\underline{g}_m}[X;K]
-{\cal{E}}_{{\tau\wedge s},\tau}^g[X;K]|^2\right]\\
&=&\lim_{m\rightarrow\infty}E\left[\sup_{s\in[0,T]}
|{\cal{E}}_{{\tau\wedge s},\tau}^{\underline{g}_m^K}[X+K_\tau]-{\cal{E}}_{{\tau\wedge s},\tau}^{{g}^K}[X+K_\tau]|^2\right]=0.\ \ \ \ \ \ \ \ \ \ \ \ \ \ \ \ \ \ \ \ \ \ \ \ \ \ \ \ \ (2.7)
\end{eqnarray*}
By (2.3), we have for each $s\in[0,T],$
\begin{eqnarray*}
\ \ \ \ \ \ \ &&{\cal{E}}_{{\tau\wedge s},\tau}^g[X_n;K^n]-{\cal{E}}_{{\tau\wedge s},\tau}^g[X;K]\\&=&{\cal{E}}_{{\tau\wedge s},\tau}^g[X_n;K^n]
-{\cal{E}}_{{\tau\wedge s},\tau}^{\overline{g}_m}[X_n;K^n]+{\cal{E}}_{{\tau\wedge s},\tau}^{\overline{g}_m}[X_n;K^n]
-{\cal{E}}_{{\tau\wedge s},\tau}^{\overline{g}_m}[X;K]\\
&&+{\cal{E}}_{{\tau\wedge s},\tau}^{\overline{g}_m}[X;K]-{\cal{E}}_{{\tau\wedge s},\tau}^g[X;K]\\
&\leq&{\cal{E}}_{{\tau\wedge s},\tau}^{\overline{g}_m}[X_n;K^n]
-{\cal{E}}_{{\tau\wedge s},\tau}^{\overline{g}_m}[X;K]+{\cal{E}}_{{\tau\wedge s},\tau}^{\overline{g}_m}[X;K]-{\cal{E}}_{{\tau\wedge s},\tau}^g[X;K],\ \ \  \ \ \  \ \ \ \ \ \ \ \ \ \ \ \ \ \ \ \ \  (2.8)
\end{eqnarray*}
and
\begin{eqnarray*}
\ \ \ \ \ \ \ &&{\cal{E}}_{{\tau\wedge s},\tau}^g[X_n;K^n]-{\cal{E}}_{{\tau\wedge s},\tau}^g[X;K]\\&=&{\cal{E}}_{{\tau\wedge s},\tau}^g[X_n;K^n]
-{\cal{E}}_{{\tau\wedge s},\tau}^{\underline{g}_m}[X_n;K^n]+{\cal{E}}_{{\tau\wedge s},\tau}^{\underline{g}_m}[X_n;K^n]
-{\cal{E}}_{{\tau\wedge s},\tau}^{\underline{g}_m}[X;K]\\
&&+{\cal{E}}_{{\tau\wedge s},\tau}^{\underline{g}_m}[X;K]-{\cal{E}}_{{\tau\wedge s},\tau}^g[X;K]\\
&\geq&{\cal{E}}_{{\tau\wedge s},\tau}^{\underline{g}_m}[X_n;K^n]
-{\cal{E}}_{{\tau\wedge s},\tau}^{\underline{g}_m}[X;K]+{\cal{E}}_{{\tau\wedge s},\tau}^{\underline{g}_m}[X;K]-{\cal{E}}_{{\tau\wedge s},\tau}^g[X;K].\ \ \ \ \ \  \ \ \ \ \ \ \ \ \ \ \ \ \ \ \ \ \  (2.9)
\end{eqnarray*}
By (2.4)-(2.9), we can complete the proof. \ \ $\Box$ \\\\
\
\textbf{Lemma 2.6} \textit{Let $g$ satisfy (A1) and (A2) with $g(s,0,0)\in{L_{\cal{F}}^\infty}(0,T)$, $K_t=\int_0^t\gamma_sds$ with $\gamma_t\in{L_{\cal{F}}^\infty}(0,T),$ $\sigma, \tau\in{\cal{T}}_{0,T}$ and $\sigma\leq\tau.$ Then for $X\in L^\infty({\cal{F}}_\tau),$ we have
$$\left\|{\cal{E}}^g_{\tau\wedge s, \tau}[X;K]\right\|_{L_{\cal{F}}^\infty(\sigma, \tau)}\leq e^{\mu\|\tau-\sigma\|_{\infty}}\left(\|X\|_{\infty}
+\|\tau-\sigma\|_{\infty}\left(\left\|g(s,0,0)\right\|_{L_{\cal{F}}^\infty(\sigma, \tau)}
+\left\|\gamma_s\right\|_{L_{\cal{F}}^\infty(\sigma, \tau)}\right)\right).$$}\\
\
\emph{Proof.} By Fan and Jiang [5, Lemma 4], we have
$$\mu|y|+\phi(|z|)\leq\mu|y|+n|z|+\phi\left(\frac{2\nu}{n}\right),\ \ \textrm{for}\ \ n\geq2\nu.\eqno(2.10)$$
Then, by (A1), we have
$$|g|\leq\mu|y|+n|z|+\phi\left(\frac{2\nu}{n}\right)+|g(s,0,0)|:=f_n(t,y,z),\ \ \textrm{ for}\ \ n\geq2\nu.\eqno(2.11)$$
For $X\in L^\infty ({\cal{F}}_\tau),$ we consider the following BSDE:
$$Y_\sigma=X+K_\tau-K_\sigma+\int_\sigma^\tau f_n(s,Y_s,Z_s)
ds-\int_\sigma^\tau Z_sdB_s,\ \ t\in[0,T].\eqno(2.12)$$
By linearization for (2.12) and $K_s=\int_0^t\gamma_sds$, we have
$$Y_\sigma=X+\int_\sigma^\tau(a_sY_s+Z_sb_s+f_n(s,0,0)+\gamma_s)
ds-\int_\sigma^\tau Z_sdB_s,\ \ t\in[0,T].\eqno(2.13)$$
where $$ a_s=\frac{f_n\left(s,Y_s,Z_s\right)-f_n\left(s,0,Z_s\right)}{Y_s}1_{|Y_s|>0}\ \ \ \textmd{and}\ \ \
b_s=\frac{(f_n\left(s,0,Z_s\right)-f_n\left(s,0,0\right))Z_s}{|Z_s|^2}1_{|Z_s|>0}.$$
Clearly, $|a_s|\leq \mu, |b_s|\leq n$ and $\|f_n(s,0,0)+\gamma_s\|_{L_{\cal{F}}^\infty(0,T)}<\infty.$

Then by the explicit solution of linear BSDE (2.13) (see Pham [17, Proposition 6.2.1]), we can get
$${\cal{E}}^{f_n}_{\sigma,\tau}[X;K]=Y_\sigma=\Gamma_\sigma^{-1}E\left[X\Gamma_\tau+\int_\sigma^\tau \Gamma_s(f_n(s,0,0)+\gamma_s)ds|{\cal{F}}_\sigma\right],\eqno(2.14)$$
where $$\Gamma_s=\exp\left\{\int_0^sb_rdB_r-\frac{1}{2}\int_0^s|b_r|^2dr+\int_0^sa_rdr\right\}.$$
Let $Q$ be a probability measure such that $\frac{dQ}{dP}=\exp\left\{\int_0^Tb_sdB_s-\frac{1}{2}\int_0^T|b_s|^2ds\right\}.$ By (2.14), we have
\begin{eqnarray*}
\left|{\cal{E}}^{f_n}_{\sigma,\tau}[X;K]\right|&=& \left\|E_Q\left[Xe^{\int_\sigma^\tau a_sds}|{\cal{F}}_\sigma\right]\right\|_{\infty}+ \left\|\int_0^T E_Q\left[1_{[\sigma,\tau]}(s)(f_n(s,0,0)+\gamma_s)e^{\int_\sigma^s a_rdr}|{\cal{F}}_\sigma\right]ds\right\|_{\infty}\\&\leq&\left\|E_Q\left[Xe^{\int_\sigma^\tau a_sds}|{\cal{F}}_\sigma\right]\right\|_{\infty}+\left\|E_Q\left[\int_\sigma^\tau (f_n(s,0,0)+\gamma_s)e^{\int_\sigma^s a_rdr}ds|{\cal{F}}_\sigma\right]\right\|_{\infty}\\
&\leq&e^{\mu\|\tau-\sigma\|_{\infty}}\left(\|X\|_{\infty}
+\|\tau-\sigma\|_{\infty}\left(\left\|f_n\left(s,0,0\right)\right\|_{L_{\cal{F}}^\infty(\sigma, \tau)}
+\left\|\gamma_s\right\|_{L_{\cal{F}}^\infty(\sigma, \tau)}\right)\right).
\end{eqnarray*}
From this, it follows that
$$\sup_{s\in[0,T]}\left|{\cal{E}}^{f_n}_{(\sigma\vee s)\wedge\tau,\tau}[X;K]\right|\leq e^{\mu\|\tau-\sigma\|_{\infty}}\left(\|X\|_{\infty}
+\|\tau-\sigma\|_{\infty}\left(\left\|f_n\left(s,0,0\right)\right\|_{L_{\cal{F}}^\infty(\sigma, \tau)}
+\left\|\gamma_s\right\|_{L_{\cal{F}}^\infty(\sigma, \tau)}\right)\right).$$
Thus we have
$$\left\|{\cal{E}}^{f_n}_{\tau\wedge s, \tau}[X;K]\right\|_{L_{\cal{F}}^\infty(\sigma, \tau)}\leq e^{\mu\|\tau-\sigma\|_{\infty}}\left(\|X\|_{\infty}
+\|\tau-\sigma\|_{\infty}\left(\left\|f_n\left(s,0,0\right)\right\|_{L_{\cal{F}}^\infty(\sigma, \tau)}
+\left\|\gamma_s\right\|_{L_{\cal{F}}^\infty(\sigma, \tau)}\right)\right).\eqno(2.15)$$
Similarly, we have
$$\left\|{\cal{E}}^{-f_n}_{\tau\wedge s, \tau}[X;K]\right\|_{L_{\cal{F}}^\infty(\sigma, \tau)}\leq e^{\mu\|\tau-\sigma\|_{\infty}}\left(\|X\|_{\infty}
+\|\tau-\sigma\|_{\infty}\left(\left\|f_n\left(s,0,0\right)\right\|_{L_{\cal{F}}^\infty(\sigma, \tau)}
+\left\|\gamma_s\right\|_{L_{\cal{F}}^\infty(\sigma, \tau)}\right)\right).\eqno(2.16)$$
On the other hand, by comparison theorem (see Jia [8, Theorem 3.6.3]), we have $\forall s\in[0,T],$
$${\cal{E}}^{-f_n}_{\tau\wedge s, \tau}[X;K]\leq{\cal{E}}^g_{\tau\wedge s, \tau}[X;K]\leq{\cal{E}}^{f_n}_{\tau\wedge s, \tau}[X;K], \  \ n\geq2\nu,\ P-a.s.\eqno(2.17)$$
Thus by (2.15)-(2.17), (2.11), the continuity of $\phi$ and $\phi(0)=0,$ we have
$$\left\|{\cal{E}}^g_{\tau\wedge s, \tau}[X;K]\right\|_{L_{\cal{F}}^\infty(\sigma, \tau)}\leq e^{\mu\|\tau-\sigma\|_{\infty}}\left(\|X\|_{\infty}
+\|\tau-\sigma\|_{\infty}\left(\left\|g(s,0,0)\right\|_{L_{\cal{F}}^\infty(\sigma, \tau)}
+\left\|\gamma_s\right\|_{L_{\cal{F}}^\infty(\sigma, \tau)}\right)\right).$$
as $n\rightarrow\infty.$ The proof is complete.\ \ $\Box$\\\\
\
\textbf{Lemma 2.7}  \textit{Let $g$ satisfies (A1) and (A2) with $g(s,0,0)\in{L_{\cal{F}}^\infty}(0,T)$, $K_s=\int_0^t\gamma_sds$ with $\gamma_t\in{L_{\cal{F}}^\infty}(0,T),$ $\sigma, \tau\in{\cal{T}}_{0,T}$ and $\sigma\leq\tau.$ Then for $X\in L^\infty({\cal{F}}_\sigma),$ we have
$$ \left\|{\cal{E}}^{g}_{\tau\wedge s,\tau}[X;K]-X\right\|_{L_{\cal{F}}^\infty(\sigma, \tau)}\leq e^{\mu\|\tau-\sigma\|_{\infty}}
\|\tau-\sigma\|_{\infty}\left(\mu\|X\|_{\infty}+\left\|g(s,0,0)\right\|_{L_{\cal{F}}^\infty(\sigma, \tau)}
+\left\|\gamma_s\right\|_{L_{\cal{F}}^\infty(\sigma, \tau)}\right).$$}\\
\
\emph{Proof.}  For $X\in L^\infty ({\cal{F}}_\sigma)$ and $s\in{[0,T]},$ set
$$g^X(s,y,z):=1_{[\sigma,\tau]}(s)g(s,y+X,z)+1_{[0,\sigma)\cup(\tau,T]}(s)g(s,y,z).\eqno(2.18)$$
Clearly, $g^X$ satisfies (A1) and (A2) with $g^X(s,0,0)\in{L_{\cal{F}}^\infty}(0,T)$.
Then by the uniqueness of solutions, we can check that for each $s\in[0,T],$
$${\cal{E}}^{g}_{(\sigma\vee s)\wedge\tau,\tau}[X;K]-X={\cal{E}}^{g^X}_{(\sigma\vee s)\wedge\tau,\tau}[0;K],\ \  P-a.s.$$
Thus by Lemma 2.6, (2.18) and (A1), we have
\begin{eqnarray*}
\left\|{\cal{E}}^{g}_{\tau\wedge s,\tau}[X;K]-X\right\|_{L_{\cal{F}}^\infty(\sigma, \tau)}
&=&\left\|{\cal{E}}^{g^X}_{\tau\wedge s,\tau}[0;K]\right\|_{L_{\cal{F}}^\infty(\sigma, \tau)}\\
&\leq& e^{\mu\|\tau-\sigma\|_{\infty}}\|\tau-\sigma\|_{\infty}\left(\|g^X(s,0,0)\|_{L_{\cal{F}}^\infty(\sigma, \tau)}
+\|\gamma_s\|_{L_{\cal{F}}^\infty(\sigma, \tau)}\right)\\
&\leq& e^{\mu\|\tau-\sigma\|_{\infty}}\|\tau-\sigma\|_{\infty}\left(\mu\|X\|_{\infty}
+\left\|g(s,0,0)\right\|_{L_{\cal{F}}^\infty(\sigma, \tau)}
+\left\|\gamma_s\right\|_{L_{\cal{F}}^\infty(\sigma, \tau)}\right).
\end{eqnarray*}
The proof is complete. $\Box$\\\\
\
\textbf{Lemma 2.8} \textit{Let $g$ satisfy (A1) and (A2) with $g(s,0,0)\in{L_{\cal{F}}^\infty}(0,T),$ $K_t=\int_0^t\gamma_sds$ with $\gamma_s\in L^\infty_{{\cal{F}}}(0,T), $ $\tau\in{\cal{T}}_{0,T}$ and $\{\tau_n\}_{n\geq1}\subset{\cal{T}}_{0,T}$ is a decreasing sequence. Let $X\in L^\infty({\cal{F}}_\tau), X_n\in L^2({\cal{F}}_{\tau_n}), n\geq1.$ If  $\|\tau_n-\tau\|_\infty\rightarrow 0$ and $X_n\rightarrow X$ in $L^2({\cal{F}}_T),$ as $n\rightarrow\infty,$ then we have
$$\lim_{n\rightarrow\infty}E\left[\sup_{s\in[0,T]}\left|{\cal{E}}_{{\tau\wedge s},\tau_n}^g\left[X_n;K\right]-{\cal{E}}_{{\tau\wedge s},\tau_n}^g\left[X;K\right]\right|^2\right]=0.$$}\\
\
\emph{Proof.}  For $m>(\mu\vee\nu),$ let $\underline{g}_m$ and $\overline{g}_m$ be defined as in (2.1) and (2.2), respectively. Firstly, we can get
\begin{eqnarray*}
\ \ \ &&\lim_{n\rightarrow\infty}\left\|\sup_{s\in[0,T]}|{\cal{E}}_{{\tau\wedge s},\tau_n}^{\underline{g}_m}[X;K]-{\cal{E}}_{{\tau\wedge s},\tau_n}^{g}[X;K]-({\cal{E}}_{{\tau\wedge s},\tau}^{\underline{g}_m}[X;K]-{\cal{E}}_{{\tau\wedge s},\tau}^{g}[X;K])|\right\|_\infty\\
&\leq&\lim_{n\rightarrow\infty}\left\|\sup_{s\in[0,T]}|{\cal{E}}_{{\tau\wedge s},\tau}^{\underline{g}_m}[{\cal{E}}_{\tau,\tau_n}^{\underline{g}_m}[X;K];K]-{\cal{E}}_{{\tau\wedge s},\tau}^{\underline{g}_m}[X;K]|\right\|_\infty\\&&+\lim_{n\rightarrow\infty}\left\|\sup_{s\in[0,T]}|{\cal{E}}_{{\tau\wedge s},\tau}^{g}[{\cal{E}}_{\tau,\tau_n}^{g}[X;K];K]-{\cal{E}}_{{\tau\wedge s},\tau}^{g}[X;K]|\right\|_\infty
\\
&\leq&\lim_{n\rightarrow\infty}\left\|\sup_{s\in[0,T]}|{\cal{E}}_{{\tau\wedge s},\tau}^{m,m}[{\cal{E}}_{\tau,\tau_n}^{\underline{g}_m}[X;K]-X]|+\sup_{s\in[0,T]}|{\cal{E}}_{{\tau\wedge s},\tau}^{-m,-m}[{\cal{E}}_{\tau,\tau_n}^{\underline{g}_m}[X;K]-X]|\right\|_\infty\\&&+\lim_{n\rightarrow\infty}\left\|\sup_{s\in[0,T]}|{\cal{E}}_{{\tau\wedge s},\tau}^{\mu,\phi}[{\cal{E}}_{\tau,\tau_n}^{g}[X;K]-X]|+\sup_{s\in[0,T]}|{\cal{E}}_{{\tau\wedge s},\tau}^{-\mu,-\phi}[{\cal{E}}_{\tau,\tau_n}^{g}[X;K]-X]|\right\|_\infty
\\
&\leq&\lim_{n\rightarrow\infty}C\|{\cal{E}}_{\tau,\tau_n}^{\underline{g}_m}[X;K]-X\|_\infty+\lim_{n\rightarrow\infty}C \|{\cal{E}}_{\tau,\tau_n}^{g}[X;K]-X\|_\infty\\
&=&0.\ \ \ \ \ \ \ \ \ \ \ \ \ \ \ \ \ \ \ \ \ \ \ \  \ \ \ \ \ \ \  \ \ \ \ \ \ \ \ \ \ \ \ \ \ \ \ \ \ \ \ \ \ \ \ \ \ \ \ \ \ \ \ \ \ \ \ \ \ \ \ \ \ \ \ \ \ \ \ \ \ \ \ \ \ \ \ \ \ \ \ \ \ \ \ \ \  \ \ \ \ \ \  \ (2.19)
\end{eqnarray*}
In the above, $C$ is a constant only dependent on $m, \mu$ and $T,$ the first inequality is due to "Consistency", the second inequality is due to the fact $\overline{g}_m$ and $\underline{g}_m$ are both Lipschitz with constant $m$ and (ii) in Remark 2.2, the third inequality is due to Lemma 2.6, the last equality is due to Lemma 2.7.

Similarly, we also have
$$\lim_{n\rightarrow\infty}\left\|\sup_{s\in[0,T]}|{\cal{E}}_{{\tau\wedge s},\tau_n}^{\overline{g}_m}[X;K]-{\cal{E}}_{{\tau\wedge s},\tau_n}^{g}[X;K]-({\cal{E}}_{{\tau\wedge s},\tau}^{\overline{g}_m}[X;K]-{\cal{E}}_{{\tau\wedge s},\tau}^{g}[X;K])|\right\|_\infty=0.\eqno(2.20)$$
Then we can complete this proof from the following inequality
\begin{eqnarray*}
&&\lim_{n\rightarrow\infty}E\left[\sup_{s\in[0,T]}|{\cal{E}}_{{\tau\wedge s},\tau_n}^g[X_n;K]-{\cal{E}}_{{\tau\wedge s},\tau_n}^g[X;K]|^2\right]\\
&\leq&\lim_{n\rightarrow\infty}2E\left[\sup_{s\in[0,T]}|{\cal{E}}_{{\tau\wedge s},\tau_n}^{\underline{g}_m}[X_n;K]
-{\cal{E}}_{{\tau\wedge s},\tau_n}^{\underline{g}_m}[X;K]+{\cal{E}}_{{\tau\wedge s},\tau_n}^{\underline{g}_m}[X;K]-{\cal{E}}_{{\tau\wedge s},\tau_n}^g[X;K]|^2\right]\\
&&+\lim_{n\rightarrow\infty}2E\left[\sup_{s\in[0,T]}|{\cal{E}}_{{\tau\wedge s},\tau_n}^{\overline{g}_m}[X_n;K]
-{\cal{E}}_{{\tau\wedge s},\tau_n}^{\overline{g}_m}[X;K]+{\cal{E}}_{{\tau\wedge s},\tau_n}^{\overline{g}_m}[X;K]-{\cal{E}}_{{\tau\wedge s},\tau_n}^g[X;K]|^2\right]
\\
&\leq&\lim_{n\rightarrow\infty}16E\left[\sup_{s\in[0,T]}|{\cal{E}}_{{\tau\wedge s},\tau_n}^{m,m}[X_n-X]|^2+\sup_{s\in[0,T]}|{\cal{E}}_{{\tau\wedge s},\tau_n}^{-m,-m}[X_n-X]|^2\right]\\&&+\lim_{n\rightarrow\infty}4E\left[\sup_{s\in[0,T]}|{\cal{E}}_{{\tau\wedge s},\tau_n}^{\underline{g}_m}[X;K]-{\cal{E}}_{{\tau\wedge s},\tau_n}^g[X;K]|^2\right]\\&&+\lim_{n\rightarrow\infty}4E\left[\sup_{s\in[0,T]}|{\cal{E}}_{{\tau\wedge s},\tau_n}^{\overline{g}_m}[X;K]-{\cal{E}}_{{\tau\wedge s},\tau_n}^g[X;K]|^2\right]
\\
&\leq&\lim_{n\rightarrow\infty}16E\left[\sup_{s\in[0,T]}|{\cal{E}}_{{\tau_n\wedge s},\tau_n}^{m,m}[X_n-X]|^2+\sup_{s\in[0,T]}|{\cal{E}}_{{\tau_n\wedge s},\tau_n}^{-m,-m}[X_n-X]|^2\right]\\&&+\lim_{m\rightarrow\infty}8E\left[\sup_{s\in[0,T]}|{\cal{E}}_{{\tau\wedge s},\tau_n}^{\underline{g}_m}[X;K]-{\cal{E}}_{{\tau\wedge s},\tau_n}^{g}[X;K]-({\cal{E}}_{{\tau\wedge s},\tau}^{\underline{g}_m}[X;K]-{\cal{E}}_{{\tau\wedge s},\tau}^{g}[X;K])|^2\right]\\&&+\lim_{m\rightarrow\infty}8E\left[\sup_{s\in[0,T]}|{\cal{E}}_{{\tau\wedge s},\tau}^{\underline{g}_m}[X;K]-{\cal{E}}_{{\tau\wedge s},\tau}^g[X;K]|^2\right]\\&&+\lim_{m\rightarrow\infty}8E\left[\sup_{s\in[0,T]}|{\cal{E}}_{{\tau\wedge s},\tau_n}^{\overline{g}_m}[X;K]-{\cal{E}}_{{\tau\wedge s},\tau_n}^{g}[X;K]-({\cal{E}}_{{\tau\wedge s},\tau}^{\overline{g}_m}[X;K]-{\cal{E}}_{{\tau\wedge s},\tau}^{g}[X;K])|^2\right]\\&&+\lim_{m\rightarrow\infty}8E\left[\sup_{s\in[0,T]}|{\cal{E}}_{{\tau\wedge s},\tau}^{\overline{g}_m}[X;K]-{\cal{E}}_{{\tau\wedge s},\tau}^g[X;K]|^2\right]\\
&\leq&\lim_{n\rightarrow\infty}32CE|X_n-X]|^2\\
&=&0.
\end{eqnarray*}
In the above, $C$ is a constant only dependent on $m$ and $T,$ the first inequality is due to the arguments of (2.8) and (2.9), the second inequality is due to the fact $\overline{g}_m$ and $\underline{g}_m$ are both Lipschitz with constant $m$ and (ii) in Remark 2.2, the third inequality is due to the fact $\tau\leq \tau_n,$  the fourth equality is due to Peng [14, Lemma 10.14, Equation (10.31)], (2.19), (2.7), (2.20) and (2.6). \ \ $\Box$
%%%%%%%%%%%%%%%%%%%%%%%%%%%%%%%%%%%%%%%%%%%%%%%%%%%%%%%%%%%%%%%%
\section{Dynamically consistent nonlinear evaluations}
%%%%%%%%%%%%%%%%%%%%%%%%%%%%%%%%%%%%%%%%%%%%%%%%%%%%%%%%%%%%%%%%
In this section, we will give the definitions of ${\cal{F}}$-evaluation $({\cal{E}}_{s,t}[\cdot])_{0\leq s\leq t\leq T}$ and related ${\cal{F}}$-evaluation $({\cal{E}}_{s,t}[\cdot,K])_{0\leq s\leq t\leq T}$ introduced by Peng [14, 16]. ${\cal{F}}$-evaluation provides an ideal characterization for the dynamical behaviors of the risk measures and the pricing of contingent claims (see Peng [14, 16] for details). \\\\
\textbf{Definition 3.1} Define a system of operators:
$${\cal{E}}_{s,t}[\cdot]:\ L^2({\cal{F}}_t)\longrightarrow L^2({\cal{F}}_s), \ 0\leq s\leq t\leq T.$$
The system is called a filtration consistent evaluation (${\cal{F}}$-evaluation for short), if it satisfies the following aximos:

(i) Monotonicity: ${\cal{E}}_{s,t}[\xi]\geq{\cal{E}}_{s,t}[\eta], P-a.s., $ if $\xi\geq \eta,\ P-a.s.;$

(ii) ${\cal{E}}_{t,t}[\xi]=\xi, P-a.s.;$

(iii) Consistency: ${\cal{E}}_{r,s}[{\cal{E}}_{s,t}[\xi]]={\cal{E}}_{r,t}[\xi], P-a.s.,$ if $r\leq s\leq t\leq T;$

(iv)  "0-1 Law": $1_A{\cal{E}}_{s,t}[\xi]=1_A{\cal{E}}_{s,t}[1_A\xi], P-a.s., $ if $A\in{\mathcal {F}}_s.$\\

Now we further give some conditions for ${\cal{F}}$-evaluation ${\cal{E}}_{s,t}[\cdot]$, where (H1) is the ${\cal{E}}_{s,t}^{\mu,\phi}$-domination property mentioned in the Introduction (see (1.2)).
\begin{itemize}
\item (H1). For each $0\leq s\leq t\leq T$ and $X,\ Y$ in $L^2({\cal{F}}_t),$ we have
\begin{center}
${\cal{E}}_{s,t}[X]-{\cal{E}}_{s,t}[Y]\leq{\cal{E}}_{s,t}^{\mu,\phi}[X-Y],\ \ P-a.s.$
\end{center}
where $\mu$ and $\phi(\cdot)$ is the constant and function given in (A1), respectively.
\item (H2). For each $0\leq s\leq t\leq T,$ we have ${\cal{E}}_{s,t}[0]=0,\ \ P-a.s.$
\end{itemize}
\textbf{Remark 3.2}  By Peng [14, Proposition 2.2], (iv) in Definition 3.1 plus (H2) is equivalent to the following (H3).
\begin{itemize}
\item (H3).  "0-1 Law": For each $0\leq s\leq t\leq T$ and $\xi\in L^2({\cal{F}}_t),$ we have
\begin{center}
$1_A{\cal{E}}_{s,t}[\xi]={\cal{E}}_{s,t}[1_A\xi],\ P-a.s.,$\  if  $A\in{\mathcal {F}}_s.$
\end{center}
\end{itemize}
\textbf{Remark 3.3} Following Peng [14, Corollary 4.4 and Proposition 4.6], we can easily check the following fact. For $K_t\in {\cal{D}}^2_{{\cal{F}}}(0,T),$ if $g$ satisfies (A1) and (A2), then ${\cal{E}}^g_{s,t}[\cdot;K]$-evaluation is an ${\cal{F}}$-evaluation and satisfies (H1). Moreover if $g$ also satisfies (A3), then we can check ${\cal{E}}^g_{s,t}[\cdot]$-evaluation satisfies (H2), thus by Remark 3.2, ${\cal{E}}^g_{s,t}[\cdot]$-evaluation further satisfies (H3).\\

Now, we give the definition of ${\cal{F}}$-expectation introduced in Coquet et al. [3] and Peng [16]. ${\cal{F}}$-expectation is a special case of ${\cal{F}}$-evaluation. For the representation of ${\cal{F}}$-expectations by solutions of BSDEs, we refer to Coquet et al. [3], Hu et al. [7] and Zheng and Li [19] for Brownian filtration and Cohen [2] and Royer [18] for general filtration.\\\\
\textbf{Definition 3.4} Define a system of operators:
$${\cal{E}}[\cdot|{\cal{F}}_t]:\ L^2({\cal{F}}_T)\longrightarrow L^2({\cal{F}}_t), \ t\in[0,T].$$
The system is called a filtration consistent condition expectation (${\cal{F}}$-expectation for short), if it satisfies the following aximos:

(i) Monotonicity: ${\cal{E}}[\xi|{\cal{F}}_t]\geq{\cal{E}}[\eta|{\cal{F}}_t], P-a.s., $ if $\xi\geq \eta,\ dP-a.s.;$

(ii) Constant preservation: ${\cal{E}}[\xi|{\cal{F}}_t]=\xi, P-a.s., $ if $\xi\in L^2({\mathcal {F}}_t);$

(iii) Consistency: ${\cal{E}}[{\cal{E}}[\xi|{\cal{F}}_t|{\cal{F}}_s]={\cal{E}}[\xi|{\cal{F}}_s], P-a.s.,$ if $s\leq t\leq T;$

(iv) "0-1 Law": ${\cal{E}}[1_A\xi|{\cal{F}}_t]=1_A{\cal{E}}[\xi|{\cal{F}}_t], P-a.s., $ if $A\in{\mathcal {F}}_t.$\\

Let ${\cal{F}}$-evaluation ${\cal{E}}_{s,t}[\cdot]$ satisfy (H1). We will introduce an ${\cal{F}}$-evaluation ${\cal{E}}_{s,t}[\cdot;K]$ generated by ${\cal{E}}_{s,t}[\cdot]$ and $K_t\in {\cal{D}}^2_{{\cal{F}}}(0,T),$ using the method in Peng [14, Section 5]. We only sketch this definition. We divide this definition into two steps.

\textbf{Step I}. Firstly, we define the space of step processes: ${\cal{D}}^{2,0}_{\cal{F}}(0,T):=\{K\in{\cal{D}}^{2}_{\cal{F}}(0,T);\ K_s=\sum_{i=0}^{N-1}\xi_i1_{[t_i,t_{i+1})}(s),$ where $t_0<t_1<\cdots<t_N$ is a partition of $[0,T]$ and $\xi_i\in L^2({\cal{F}}_{t_i})\}.$ As Peng [14, Definition 5.2 and Lemma 5.4],  we have the following Proposition 3.5.\\\\
\textbf{Proposition 3.5} \textit{Let ${\cal{F}}$-evaluation ${\cal{E}}_{s,t}[\cdot]$ satisfy (H1). For each $K_t\in {\cal{D}}^{2,0}_{{\cal{F}}}(0,T)$ with form $K_s=\sum_{i=0}^{N-1}\xi_i1_{[t_i,t_{i+1})}(s)$, where $t_0<t_1<\cdots<t_N$ is a partition of $[0,T]$ and $\xi_i\in L^2({\cal{F}}_{t_i})$, there exists a unique ${\cal{F}}$-evaluation, denoted by ${\cal{E}}_{s,t}[\cdot;K]$ such that $\forall t_i\leq s\leq t\leq t_{i+1}$ and $X\in L^2({\cal{F}}_t),$
$${\cal{E}}_{s,t}[X;K]={\cal{E}}_{s,t}[X+K_t-K_s],\ \ P-a.s.\eqno(3.1)$$
and for each $K,\ K'\in {\cal{D}}^{2,0}_{{\cal{F}}}(0,T)$ and $0\leq s \leq t\leq T,\ X,\ X'\in L^{2}({\cal{F}}_t),$ we have
$${\cal{E}}^{-\mu,-\phi}_{s,t}[X-X';K-K']\leq{\cal{E}}_{s,t}[X;K]-{\cal{E}}_{s,t}[X';K']\leq{\cal{E}}_{s,t}^{\mu,\phi}[X-X';K-K'],\ \ P-a.s.$$}

We further have the following consequence.\\\\
\textbf{Proposition  3.6} \textit{Let ${\cal{F}}$-evaluation ${\cal{E}}_{s,t}[\cdot]$ satisfy (H1) and $K^n\in {\cal{D}}^{2,0}_{{\cal{F}}}(0,T), n\geq1, t\in [0,T].$ If $\{K^n\}_{n\geq1}$ is a Cauchy sequence in $L^2_{\cal{F}}(0,t),$ $\{K^n_t\}_{n\geq1}$ and $\{X_n\}_{n\geq1}$ are both Cauchy sequences in $L^2({\cal{F}}_T),$ then we have
$$\lim_{m, n\rightarrow\infty}E\left[\sup_{0\leq s\leq t}|{\cal{E}}_{s,t}[X_m;K^m]+K_s^m-{\cal{E}}_{s,t}[X_n;K^n]-K_s^n|^2\right]=0.$$}
\emph{Proof.} By Proposition 3.5, Lemma 2.5 and the fact ${\cal{E}}_{s,t}^{\mu,\phi}[0;0]={\cal{E}}^{-\mu,-\phi}_{s,t}[0;0]=0$, we have
\begin{eqnarray*}
&&\lim_{m, n\rightarrow\infty}E\left[\sup_{0\leq s\leq t}|{\cal{E}}_{s,t}[X_m;K^m]+K_s^m-{\cal{E}}_{s,t}[X_n;K^n]-K_s^n|^2\right]\\
&\leq&\lim_{m, n\rightarrow\infty}2E\left[\sup_{0\leq s\leq t}|{\cal{E}}_{s,t}^{\mu,\phi}[X_m-X_n;K^m-K^n]+K_s^m-K_s^n|^2\right]\\
&&+\lim_{m, n\rightarrow\infty}2E\left[\sup_{0\leq s\leq t}|{\cal{E}}^{-\mu,-\phi}_{s,t}[X_m-X_n;K^m-K^n]+K_s^m-K_s^n|^2\right]\\
&=&0.
\end{eqnarray*}
The proof is complete.\ \ $\Box$ \\

\textbf{Step II}. For $K\in {\cal{D}}^{2}_{\cal{F}}(0,T)$ and $\forall 0\leq s\leq t\leq T,$ by Peng [14, Remark 5.5.1], we can take partitions $0=t_0^i<t_1^i<\cdots<t_i^i=T$ of $[0,T], i\geq1$ such that $\max_j(t^i_{j+1}-t^i_j)\rightarrow0$  with $s=t^i_{j_1}$ and $t=t^i_{j_2},$ for some $j_1\leq j_2\leq i.$ We define $K_s^i:=\sum_{j=0}^{i-1}K_{t_j^i}1_{[t^i_j,t^i_{j+1})}(s).$
Thus $K^i$ converges to $K$ in $L^2_{\cal{F}}(0,T)$ and  $K_s^i=K_s,$ $K_t^i=K_t.$  Then for $X\in L^2({\cal{F}}_t),$ by Proposition 3.6, we can get $\{{\cal{E}}_{s,t}[X;K^i]\}_{i\geq1}$ is a Cauchy sequence in $L^2({\cal{F}}_T).$ We define
$${\cal{E}}_{s,t}[X;K]:=\lim_{i\rightarrow\infty}{\cal{E}}_{s,t}[X;K^i]\ \ \textrm{in} \ \ L^2({\cal{F}}_T).$$
The Definition of ${\cal{E}}_{s,t}[\cdot;K]$ is complete.\\

By Definition of ${\cal{E}}_{s,t}[\cdot;K]$, Proposition 3.5 and Lemma 2.5, we can get Proposition 3.7, immediately. We omit its proof.\\\\
\textbf{Proposition 3.7} \textit{Let ${\cal{F}}$-evaluation ${\cal{E}}_{s,t}[\cdot]$ satisfy (H1). Then for each $K_t\in {\cal{D}}^{2}_{{\cal{F}}}(0,T),$ ${\cal{E}}_{s,t}[\cdot;K]$ is an ${\cal{F}}$-evaluation, such that for $K,\ K'\in {\cal{D}}^{2}_{{\cal{F}}}(0,T),$ $t\in[0,T]$ and $X,\ X'\in L^2({\cal{F}}_t),$ we have for $s\in[0,t],$ $P-a.s.,$
$${\cal{E}}^{-\mu,-\phi}_{s,t}[X-X';K-K']\leq{\cal{E}}_{s,t}[X;K]-{\cal{E}}_{s,t}[X';K']
\leq{\cal{E}}_{s,t}^{\mu,\phi}[X-X';K-K'],\eqno(3.2)$$}

For ${\cal{F}}$-evaluation ${\cal{E}}_{s,t}[\cdot;K],$ we further have the the following properties.\\\\
\
\textbf{Corollary 3.8}\textit{ Let ${\cal{F}}$-evaluation ${\cal{E}}_{s,t}[\cdot]$ satisfy (H1) and (H2), $K_t,\ K'_t\in {\cal{D}}^2_{{\cal{F}}}(0,T).$ Then for each $t\in [0,T]$ and $X$ in $L^2({\cal{F}}_t),$ we have $\forall s\in[0,t],$}

\textit{(i) ${\cal{E}}_{s,t}^{-\mu,-\phi}[X;K]\leq{\cal{E}}_{s,t}[X;K]\leq{\cal{E}}_{s,t}^{\mu,\phi}[X;K],\ \ P-a.s.;$
}

\textit{(ii) $|{\cal{E}}_{s,t}[X]|\leq{\cal{E}}_{s,t}^{\mu,\phi}[|X|],\ \ P-a.s.$}\\\\
\
\emph{Proof.} By (3.1), we have $\forall s\in[0,t],$
$${\cal{E}}_{s,t}[X;0]={\cal{E}}_{s,t}[X],\ \ P-a.s.\eqno(3.3)$$
By (3.3), (H2) and (3.2), we have $\forall s\in[0,t],$ $P-a.s.,$
$${\cal{E}}^{-\mu,-\phi}_{s,t}[X;K]={\cal{E}}^{-\mu,-\phi}_{s,t}[X;K]+{\cal{E}}_{s,t}[0;0]
\leq{\cal{E}}_{s,t}[X;K]
\leq{\cal{E}}_{s,t}^{\mu,\phi}[X;K]+{\cal{E}}_{s,t}[0;0]={\cal{E}}^{\mu,\phi}_{s,t}[X;K].$$
Then we obtain (i). We can easily check $\forall s\in[0,t],$
$$-{\cal{E}}_{s,t}^{\mu,\phi}[X;K]={\cal{E}}_{s,t}^{-\mu,-\phi}[-X;-K],\ \ P-a.s.$$
By this, comparison theorem (Jia [9, Theorem 3.1]), (i) and (3.3), we have $\forall s\in[0,t],$
$$-{\cal{E}}_{s,t}^{\mu,\phi}[|X|]={\cal{E}}_{s,t}^{-\mu,-\phi}[-|X|]\leq{\cal{E}}_{s,t}^{-\mu,-\phi}[X]\leq{\cal{E}}_{s,t}[X]
\leq{\cal{E}}_{s,t}^{\mu,\phi}[X]\leq{\cal{E}}_{s,t}^{\mu,\phi}[|X|],\ \ P-a.s.$$
Thus, (ii) is true. The proof is complete. $\Box$\\\\
\textbf{Lemma 3.9} \textit{Let ${\cal{F}}$-evaluation ${\cal{E}}_{s,t}[\cdot]$ satisfy (H1), $K_t,\ K^n_t\in {\cal{D}}^2_{{\cal{F}}}(0,T),$ $t\in [0,T]$ and $X,\ X_n$ in $L^2({\cal{F}}_t),\ n\geq1,$
If $K^n\rightarrow K$ in $L^2_{\cal{F}}(0,T),$ $K^n_t\rightarrow K_t$ and $X_n\rightarrow X$ both in $L^2({\cal{F}}_T),$ as $n\rightarrow\infty,$ then we have}
$$\lim_{n\rightarrow\infty}E\left[\sup_{0\leq s\leq t}|{\cal{E}}_{s,t}[X;K]+K_s-{\cal{E}}_{s,t}[X_n;K^n]-K_s^n|^2\right]=0.$$
\emph{Proof.} By (3.2) and the proof of Proposition 3.6, we can complete this proof.\ \ $\Box$ \\\\
\
\textbf{Definition 3.10} Let $K_t\in {\cal{D}}^2_{{\cal{F}}}(0,T).$ A process $Y_t$ with $Y_t\in L^2({\cal{F}}_t)$ for $t\in [0,T],$ is called an ${\cal{E}}_{s,t}[\cdot;K]$-martingale (resp. ${\cal{E}}_{s,t}[\cdot;K]$-supermartingale, ${\cal{E}}_{s,t}[\cdot;K]$-submartingale), if, for each $0\leq s\leq t\leq T,$ we have
\begin{center}
${\cal{E}}_{s,t}[Y_t;K]=Y_s,$\ \ \ (resp. $\leq,\ \geq$).
\end{center}\ \\
\textbf{Lemma 3.11}\textit{ Let ${\cal{F}}$-evaluation ${\cal{E}}_{s,t}[\cdot]$ satisfy (H1) and (H2). Then for each $t\in[0,T]$ and $X\in L^2({\cal{F}}_t),$ ${\cal{E}}_{s,t}[X]$ admits a RCLL version.}
\\\\
\
\emph{Proof.} Given $t\in[0,T]$.  As (2.1) and (2.2), we can find two functions $g_i(y,z): {\mathbf{R}}\times{\mathbf{R}}^{\mathit{d}}\mapsto {\mathbf{R}}, i=1,2,$ which both satisfy (A2) and are both Lipschitz in $(y, z)$ with some constant $C_0,$ such that for each $(y,z)\in{\mathbf{R}}\times{\mathbf{R}}^{\mathit{d}}$,
$${g}_1\leq {-\mu|y|-\phi(|z|)}\ \ \textmd{and}\ \ {g}_2\geq {\mu|y|+\phi(|z|)}.$$
By (i) in Corollary 3.8 and comparison theorem (see Jia [9, Theorem 3.1]), we have for each $X\in L^2({\cal{F}}_t)$ and $s\in[0,t]$ $${\cal{E}}_{s,t}^{g_2}[X]\geq{\cal{E}}_{s,t}^{\mu,\phi}[X]\geq {\cal{E}}_{s,t}[X]\geq{\cal{E}}_{s,t}^{-\mu,-\phi}[X]\geq {\cal{E}}_{s,t}^{g_1}[X],\ \ P-a.s.\eqno(3.4)$$
Then we can check that ${\cal{E}}_{s,t}[X]$ is an ${\cal{E}}_{s,t}^{g_1}[\cdot]$-supermartingale. Thus, by Peng [16, Theorem 3.7], we get that for a denumerable dense subset ${{D}}$ of $[0,t]$, almost all $\omega\in\Omega$ and all $r\in[0,t]$, we have $\lim_{s\in{{D}},\ s\searrow r}{\cal{E}}_{s,t}[X]$ and $\lim_{s\in{{D}},\ s\nearrow r}{\cal{E}}_{s,t}[X]$ both exist and are finite. For each $r\in[0,t)$, we set
$$Y_r:=\lim_{s\in{{D}},\ s\searrow r}{\cal{E}}_{s,t}[X],\eqno(3.5)$$
then from some classic arguments, $Y_r$ is RCLL. Thus we only need prove ${\cal{E}}_{r,t}[X]=Y_r, P-a.s.$ for $r\in[0,t).$ By (ii) in Corollary 3.8 and Jia [9, Theorem 2.3], we have
$$ E\left[\sup_{0\leq s\leq t}|{\cal{E}}_{s,t}[X]|^2\right]\leq E\left[\sup_{0\leq s\leq t}|{\cal{E}}_{s,t}^{\mu,\phi}[|X|]|^2\right]<+\infty.\eqno(3.6)$$
By (3.5), (3.6) and Lebesgue dominated convergence theorem, we have
$$\lim_{s\in{{D}},\ s\searrow r}{\cal{E}}_{s,t}[X]=Y_r, \ \ \ r\in[0,t). \eqno(3.7)$$
in $L^2({\cal{F}}_T)$ sense. By (3.4) and Peng [16, Lemma 7.6], we have
$$\lim_{s\in{{D}},\ s\searrow r}E\left[|{\cal{E}}_{r,s}[Y_r]-Y_r|^2\right]\leq\lim_{s\in{{D}},\ s\searrow r}2E\left[|{\cal{E}}_{r,s}^{g_1}[Y_r]-Y_r|^2\right]+\lim_{s\in{{D}},\ s\searrow r}2E\left[|{\cal{E}}_{r,s}^{g_2}[Y_r]-Y_r|^2\right]=0.\eqno(3.8)$$
We also have for $r\in[0,t),$
\begin{eqnarray*}
\ \ &&\lim_{s\in{{D}},\ s\searrow r}E\left[|{\cal{E}}_{r,s}^{g_2}\left[\left|{\cal{E}}_{s,t}[X]-Y_r|\right]\right|^2\right]\\
&\leq&\lim_{s\in{{D}},\ s\searrow r}CE\left[|{\cal{E}}_{s,t}[X]-Y_r|^2+\left(\int_r^s|g_2(u,0,0)|du\right)^2\right]\\
&\leq&\lim_{s\in{{D}},\ s\searrow r}CE\left[(s-r)\left(\int_r^s|g_2(u,0,0)|^2du\right)\right]\\
&=&0,\ \ \ \ \ \ \ \ \ \ \ \ \ \ \ \ \ \ \ \ \ \ \ \ \ \ \  \ \ \ \ \ \ \ \ \  \ \ \ \ \ \ \ \  \ \ \ \ \ \ \ \ \ \ \ \  \ \ \ \ \ \ \ \ \ \ \ \ \ \ \ \ \ \ \ \ \ \ \ \ \ \ \ \ \ \ \ \ \ \ \ \ \ \ \ \ \ \ \ \ \ \ \ \ \ (3.9)
\end{eqnarray*}
where $C$ is a constant only dependent on $T$ and $C_0$. In (3.9), the first inequality is from an element estimate of BSDE (see Briand et al. [1, Proposition 2.2]),  the second inequality is from (3.7) and Cauchy-Schwarz inequality, the equality is due to the fact $g_2$ satisfies (A2).

By "Consistency" of ${\cal{E}}_{r,t}[\cdot]$, (ii) in Corollary 3.8 and (3.4), we have $P-a.s.,$
\begin{eqnarray*}
\ \ \ \ \ \ \ \ |{\cal{E}}_{r,t}[X]-Y_r|&=&|{\cal{E}}_{r,s}[{\cal{E}}_{s,t}[X]]-Y_r|\\&=&|{\cal{E}}_{r,s}[{\cal{E}}_{s,t}[X]]-{\cal{E}}_{r,s}[Y_r]
+{\cal{E}}_{r,s}[Y_r]-Y_r|\\
&\leq&|{\cal{E}}_{r,s}[{\cal{E}}_{s,t}[X]]-{\cal{E}}_{r,s}[Y_r]|
+|{\cal{E}}_{r,s}[Y_r]-Y_r|\\
&\leq&{\cal{E}}_{r,s}^{\mu,\phi}[|{\cal{E}}_{s,t}[X]-Y_r|]+|{\cal{E}}_{r,s}[Y_r]-Y_r|\\
&\leq&{\cal{E}}_{r,s}^{g_2}[|{\cal{E}}_{s,t}[X]-Y_r|]+|{\cal{E}}_{r,s}[Y_r]-Y_r|. \ \ \ \ \ \ \ \ \ \ \ \ \ \ \ \ \ \ \ \ \ \ \ \ \ \ \ \ \ \ \ \ \ \ \ \ \ \ \ \ (3.10)
\end{eqnarray*}
By (3.8)-(3.10), we get that for $r\in[0,t),$ ${\cal{E}}_{r,t}[X]=Y_r,\ P-a.s.$ The proof is complete.  $\Box$\\

We will always take a RCLL version of ${\cal{E}}_{r,t}[\cdot].$ Furthermore, we have\\\\
\textbf{Corollary 3.12} \textit{Let ${\cal{F}}$-evaluation ${\cal{E}}_{s,t}[\cdot]$ satisfy (H1), (H2) and $K\in {\cal{D}}^2_{{\cal{F}}}(0,T).$  Then for each $t\in[0,T]$ and $X\in L^2({\cal{F}}_t),$ ${\cal{E}}_{s,t}[X;K]\in {\cal{D}}^2_{{\cal{F}}}(0,t).$}
\\\\
\emph{Proof. } For $K\in {\cal{D}}^{2,0}_{{\cal{F}}}(0,T),$ by (3.1), "Consistency" and Lemma 3.11, we can prove ${\cal{E}}_{s,t}[X;K]$ is RCLL. By this and Lemma 3.9, for $K\in {\cal{D}}^{2}_{{\cal{F}}}(0,T),$ we can get ${\cal{E}}_{s,t}[X;K]+K_s$ is RCLL. Thus ${\cal{E}}_{s,t}[X;K]$ is RCLL. In view of (i) in Corollary 3.8, we have ${\cal{E}}_{s,t}[X;K]\in {\cal{D}}^2_{{\cal{F}}}(0,t).$ \ \  $\Box$

%%%%%%%%%%%%%%%%%%%%%%%%%%%%%%%%%%%%%%%%%%%%%%%%%%%%%%%%%%%%%%%%
\section{Optional stopping theorem of ${\cal{E}}_{s,t}[\cdot]$-supermartingales}
%%%%%%%%%%%%%%%%%%%%%%%%%%%%%%%%%%%%%%%%%%%%%%%%%%%%%%%
In this section, we will firstly extend the definition of ${\cal{F}}$-evaluation ${\cal{E}}_{s,t}[\cdot]$ to ${\cal{E}}_{\sigma,\tau}[\cdot]$ with $\sigma, \tau\in{\cal{T}}_{0,T}.$ We divide this extension into three steps.

\textbf{Step I.}  Let ${\cal{F}}$-evaluation ${\cal{E}}_{s,t}[\cdot]$ satisfy (H1) and (H2). By the same argument as Peng [14, Section 10], we can firstly extend the definition of ${\cal{F}}$-evaluation ${\cal{E}}_{s,t}[\cdot]$ and ${\cal{E}}_{s,t}[\cdot;K]$ to ${\cal{E}}_{\sigma,\tau}[\cdot]$ and ${\cal{E}}_{\sigma,\tau}[\cdot;K]$ with $\sigma\in{\cal{T}}_{0,T}$ and $\tau\in{\cal{T}}^0_{0,T}$ for $L^2$ terminal variable. Similarly, we can obtain the following result as Peng [14, Lemma 10.13]. \\\\
\textbf{Lemma 4.1} \textit{ The system of operators
$${\cal{E}}_{\sigma,\tau}[\cdot]:\ L^2({\cal{F}}_\tau)\longrightarrow L^2({\cal{F}}_\sigma),\ \ \sigma\leq\tau,\ \sigma\in{\cal{T}}_{0,T},\ \tau\in{\cal{T}}^0_{0,T},$$
satisfy}

\textit{(i) Monotonicity: ${\cal{E}}_{\sigma,\tau}[\xi]\geq{\cal{E}}_{\sigma,\tau}[\eta], P-a.s., $ if $\xi, \eta \in L^2({\cal{F}}_\tau)$ and $\xi\geq \eta,\ P-a.s.;$}

\textit{(ii)  ${\cal{E}}_{\tau,\tau}[\xi]=\xi, P-a.s., $ if $\xi\in L^2({\mathcal {F}}_\tau);$}

\textit{(iii) Consistency: ${\cal{E}}_{\sigma,\rho}[{\cal{E}}_{\rho,\tau}[\xi]]={\cal{E}}_{\sigma,\tau}[\xi], P-a.s.,$ if $\sigma\leq \rho\leq \tau$ and $\xi\in L^2({\mathcal {F}}_\tau), \ \rho\in{\cal{T}}^0_{0,T};$}

\textit{(iv) "0-1 Law":  $1_A{\cal{E}}_{\sigma,\tau}[\xi]={\cal{E}}_{\sigma,\tau}[1_A\xi], P-a.s., $ if $A\in{\mathcal {F}}_\sigma,\ \xi\in L^2({\mathcal {F}}_\tau);$}

\textit{(v) For $K\in {\cal{D}}^2_{{\cal{F}}}(0,T),$ ${\cal{E}}_{\sigma,\tau}[\cdot;K]$ satisfies the above (i)-(iii) with ${\cal{E}}_{\sigma,\tau}[\cdot;0]={\cal{E}}_{\sigma,\tau}[\cdot]$ and}
$$1_A{\cal{E}}_{\sigma,\tau}[\xi;K]=1_A{\cal{E}}_{\sigma,\tau}[1_A\xi;K],\ \ P-a.s. \ \ \textrm{if} \ \ A\in{\mathcal {F}}_\sigma,\ \xi\in L^2({\cal{F}}_\tau);\eqno(4.1)$$

\textit{(vi) For $K\in {\cal{D}}^2_{{\cal{F}}}(0,T)$ and $\xi\in L^2({\mathcal {F}}_\tau),$ ${\cal{E}}_{\tau\wedge \cdot,\tau}[\xi;K]$ is RCLL and for $X, X'\in L^2({\mathcal {F}}_\tau)$ and $K, K'\in{\cal{D}}^2_{{\cal{F}}}(0,T),$ we have}
$${\cal{E}}^{-\mu,-\phi}_{\sigma,\tau}[X-X';K-K']\leq{\cal{E}}_{\sigma,\tau}[X;K]-{\cal{E}}_{\sigma,\tau}[X';K']
\leq{\cal{E}}_{\sigma,\tau}^{\mu,\phi}[X-X';K-K'],\ \ P-a.s.\eqno(4.2)$$

\textbf{Step II.}  In this step, we will extend the definition of ${\cal{F}}$-evaluation ${\cal{E}}_{s,t}[\cdot]$ to ${\cal{E}}_{\sigma,\tau}[\cdot],$ with $\sigma, \tau\in{\cal{T}}_{0,T}$ for bounded terminal variable. We need the following convergence results. \\\\
\
\textbf{Lemma 4.2}  \textit{Let ${\cal{F}}$-evaluation ${\cal{E}}_{s,t}[\cdot]$ satisfy (H1) and (H2). Let $\tau\in{\cal{T}}_{0,T}$ and $\{\tau_n\}_{n\geq1}\subset{\cal{T}}^0_{0,T}$ be a decreasing sequence such that for each $n\geq1,$ $\tau_n\geq\tau.$ Then we have}

\textit{(i) If $K\in {\cal{D}}^2_{{\cal{F}}}(0,T),$ $X\in L^\infty ({\cal{F}}_{\tau}),$ $X_n\in L^\infty ({\cal{F}}_{\tau_n}), n\geq1,$ and $X_n\rightarrow X$ in $L^\infty({\cal{F}}_T),$ as $n\rightarrow\infty,$ then we have
$$\lim_{n\rightarrow\infty}\left\|\sup_{t\in[0,T]}|{\cal{E}}_{\tau_n\wedge t,\tau_n}[X_n;K]-{\cal{E}}_{\tau_n\wedge t,\tau_n}[X;K]|\right\|_\infty=0.$$}

\textit{(ii) If $K\in {\cal{D}}^2_{{\cal{F}}}(0,T),$ $X\in L^2 ({\cal{F}}_{\tau}),$ $X_n\in L^2 ({\cal{F}}_{\tau_n}), n\geq1,$ and $X_n\rightarrow X$ in $L^2({\cal{F}}_T)$ and $\|\tau_n-\tau\|_\infty\rightarrow 0,$ as $n\rightarrow\infty,$ then we have}
$$\lim_{n\rightarrow\infty}E\left[\sup_{t\in[0,T]}\left|{\cal{E}}_{{\tau\wedge t},\tau_n}\left[X_n;K\right]-{\cal{E}}_{{\tau\wedge t},\tau_n}\left[X;K\right]\right|^2\right]=0.$$

\textit{(iii) If $K_t=\int_0^t\gamma_sds$ with $\gamma_s\in L^\infty_{{\cal{F}}}(0,T),$ $X\in L^\infty ({\cal{F}}_{\tau}),$ and $\|\tau_n-\tau\|_\infty\rightarrow 0,$ as $n\rightarrow\infty,$ then we have}
$$ \lim_{m,n\rightarrow\infty}\left\|\sup_{t\in[0,T]}|{\cal{E}}_{\tau\wedge t,\tau_n}[X;K]
-{\cal{E}}_{\tau\wedge t,\tau_m}[X;K]|\right\|_{\infty}=0.$$
\emph{Proof.} By (4.2), we have
\begin{eqnarray*}
&&\left\|\sup_{t\in[0,T]}|{\cal{E}}_{\tau_n\wedge t,\tau_n}[X_n;K]-{\cal{E}}_{\tau_n\wedge t,\tau_n}[X;K]\right\|_{\infty}\\
&\leq&
\left\|\sup_{t\in[0,T]}|{\cal{E}}_{\tau_n\wedge t,\tau_n}^{\mu,\phi}[X_n-X]|\right\|_{\infty}
+\left\|\sup_{t\in[0,T]}|{\cal{E}}_{\tau_n\wedge t,\tau_n}^{-\mu,-\phi}[X_n-X]|\right\|_{\infty}.
\end{eqnarray*}
Then by Lemma 2.6, we obtain (i). By (4.2), we have
\begin{eqnarray*}
&&E\left[\sup_{t\in[0,T]}\left|{\cal{E}}_{{\tau\wedge t},\tau_n}\left[X_n;K\right]-{\cal{E}}_{{\tau\wedge t},\tau_n}\left[X;K\right]\right|^2\right]\\
&\leq&2E\left[\sup_{t\in[0,T]}\left|{\cal{E}}_{{\tau\wedge t},\tau_n}^{\mu,\phi}\left[X_n-X\right]\right|^2\right]+2E\left[\sup_{t\in[0,T]}\left|{\cal{E}}_{{\tau\wedge t},\tau_n}^{-\mu,-\phi}\left[X_n-X\right]\right|^2\right].
\end{eqnarray*}
Then by Lemma 2.8, we obtain (ii). By "Consistency",  (4.2) and Lemma 2.6, we can deduce
\begin{eqnarray*}
&&\left\|\sup_{t\in[0,T]}|{\cal{E}}_{\tau\wedge t,\tau_m}[X;K]-{\cal{E}}_{\tau\wedge t,\tau_n}[X;K]|\right\|_{\infty}\\&\leq&
\left\|\sup_{t\in[0,T]}|{\cal{E}}_{\tau\wedge t,\tau_m\wedge\tau_n}[{\cal{E}}_{\tau_m\wedge\tau_n,\tau_m}[X;K];K]
-{\cal{E}}_{\tau\wedge t,\tau_m\wedge\tau_n}[X;K]|\right\|_{\infty}\\&&+
\left\|\sup_{t\in[0,T]}|{\cal{E}}_{\tau\wedge t,\tau_m\wedge\tau_n}[X;K]-{\cal{E}}_{\tau\wedge t,\tau_m\wedge\tau_n}
[{\cal{E}}_{\tau_m\wedge\tau_n,\tau_n}[X;K];K]|
\right\|_{\infty}\\
&\leq&\left\|\sup_{t\in[0,T]}|{\cal{E}}^{\mu,\phi}_{\tau\wedge t,\tau_m\wedge\tau_n}[
{\cal{E}}_{\tau_m\wedge\tau_n,\tau_m}[X;K]-X]|\right\|_{\infty}
+\left\|\sup_{t\in[0,T]}|{\cal{E}}^{-\mu,-\phi}_{\tau\wedge t,\tau_m\wedge\tau_n}[
{\cal{E}}_{\tau_m\wedge\tau_n,\tau_m}[X;K]-X]|\right\|_{\infty}\\
&&+\left\|\sup_{t\in[0,T]}|{\cal{E}}^{\mu,\phi}_{\tau\wedge t,\tau_m\wedge\tau_n}[
{\cal{E}}_{\tau_m\wedge\tau_n,\tau_n}[X;K]-X]|\right\|_{\infty}
+\left\|\sup_{t\in[0,T]}|{\cal{E}}^{-\mu,-\phi}_{\tau\wedge t,\tau_m\wedge\tau_n}[
{\cal{E}}_{\tau_m\wedge\tau_n,\tau_n}[X;K]-X]|\right\|_{\infty}\\
&\leq&2e^{\mu T}\left(\|{\cal{E}}_{\tau_m\wedge\tau_n,\tau_m}[X;K]-X\|_{\infty}+\|{\cal{E}}_{\tau_m\wedge\tau_n,\tau_n}[X;K]
-X\|_{\infty}\right)\\
&\leq&2e^{\mu T}(\|{\cal{E}}^{\mu,\phi}_{\tau_m\wedge\tau_n,\tau_m}[X;K]-X\|_{\infty}+
\|{\cal{E}}^{-\mu,-\phi}_{\tau_m\wedge\tau_n,\tau_m}[X;K]-X\|_{\infty}\\
&&+2e^{\mu T}(\|{\cal{E}}^{\mu,\phi}_{\tau_m\wedge\tau_n,\tau_n}[X;K]-X\|_{\infty}+
\|{\cal{E}}^{-\mu,-\phi}_{\tau_m\wedge\tau_n,\tau_n}[X;K]-X\|_{\infty}).
\end{eqnarray*}
Then by Lemma 2.7, we can obtain (iii).  The proof is complete. \ \ $\Box$\\

By (iii) in Lemma 4.2, the following Definition 4.3 is well defined.\\\\
\textbf{Definition 4.3 } Let ${\cal{F}}$-evaluation ${\cal{E}}_{s,t}[\cdot]$ satisfy (H1) and (H2), $K_t=\int_0^t\gamma_sds$ with $\gamma_s\in L^\infty_{{\cal{F}}}(0,T).$ Let $\sigma, \tau\in{\cal{T}}_{0,T}, \sigma\leq\tau$ and $\{\tau_n\}_{n\geq1}\subset{\cal{T}}^0_{0,T}$ be a decreasing sequence such that $\|\tau_n-\tau\|_\infty\rightarrow 0,$ as $n\rightarrow\infty.$ If $X\in L^\infty ({\cal{F}}_{\tau}),$ then we define
$${\cal{E}}_{\sigma,\tau}[X;K]:=\lim_{n\rightarrow\infty}{\cal{E}}_{\sigma,\tau_n}[X;K]\ \ \textrm{in}\ \ L^\infty({\cal{F}}_T),$$
and $${\cal{E}}_{\sigma,\tau}[X]:={\cal{E}}_{\sigma,\tau}[X;0].$$
\textbf{Lemma 4.4}  \textit{The system of operators
$${\cal{E}}_{\sigma,\tau}[\cdot]:\ L^\infty({\cal{F}}_\tau)\longrightarrow L^\infty({\cal{F}}_\sigma),\ \ \sigma\leq\tau,\ \sigma, \tau\in{\cal{T}}_{0,T},$$
satisfy}

\textit{(i) Monotonicity: ${\cal{E}}_{\sigma,\tau}[\xi]\geq{\cal{E}}_{\sigma,\tau}[\eta], P-a.s., $ if $\xi, \eta \in L^\infty({\cal{F}}_\tau)$ and $\xi\geq \eta,\ P-a.s.;$}

\textit{(ii)  ${\cal{E}}_{\tau,\tau}[\xi]=\xi, P-a.s., $ if $\xi\in L^\infty({\mathcal {F}}_\tau);$}

\textit{(iii) Consistency: ${\cal{E}}_{\sigma,\rho}[{\cal{E}}_{\rho,\tau}[\xi]]={\cal{E}}_{\sigma,\tau}[\xi], P-a.s.,$ if $\sigma\leq \rho\leq \tau$ and $\xi\in L^\infty({\mathcal {F}}_\tau), \ \rho\in{\cal{T}}_{0,T};$}

\textit{(iv) "0-1 Law":  $1_A{\cal{E}}_{\sigma,\tau}[\xi]={\cal{E}}_{\sigma,\tau}[1_A\xi], P-a.s., $ if $A\in{\mathcal {F}}_\sigma,\ \xi\in L^\infty({\mathcal {F}}_\tau);$}

\textit{(v) For $K_t=\int_0^t\gamma_sds$ with $\gamma_s\in L^\infty_{{\cal{F}}}(0,T),$ ${\cal{E}}_{\sigma,\tau}[\cdot;K]$ satisfies the above (i)-(iii) and}
$$1_A{\cal{E}}_{\sigma,\tau}[\xi;K]=1_A{\cal{E}}_{\sigma,\tau}[1_A\xi;K],\ \ P-a.s. \ \ \textrm{if} \ \ A\in{\mathcal {F}}_\sigma,\ \xi\in L^\infty({\mathcal {F}}_\tau);$$

\textit{(vi) For $K_t=\int_0^t\gamma_sds$ with $\gamma_s\in L^\infty_{{\cal{F}}}(0,T)$ and $\xi\in L^\infty({\mathcal {F}}_\tau),$ ${\cal{E}}_{\tau\wedge \cdot,\tau}[\xi;K]$ is RCLL and for $X, X'\in L^\infty({\mathcal {F}}_\tau)$ and $K'_t=\int_0^t\gamma'_sds$ with $\gamma'_s\in L^\infty_{{\cal{F}}}(0,T),$  we have
$${\cal{E}}^{-\mu,-\phi}_{\sigma,\tau}[X-X';K-K']\leq{\cal{E}}_{\sigma,\tau}[X;K]-{\cal{E}}_{\sigma,\tau}[X';K']
\leq{\cal{E}}_{\sigma,\tau}^{\mu,\phi}[X-X';K-K'],\ \ P-a.s.\eqno(4.3)$$}
\emph{Proof.} For $\tau\in{\cal{T}}_{0,T},$ we can find a decreasing sequence $\{\tau_n\}_{n\geq1}\subset{\cal{T}}^0_{0,T}$ such that $\|\tau_n-\tau\|_\infty\rightarrow 0,$ as $n\rightarrow\infty,$ by setting
$$\tau_n:=T2^{-n}i\sum_{i=1}^{2^{n}}1_{\{2^{-n}(i-1)T\leq\tau< 2^{-n}{i}T\}}+1_{\{\tau=T\}}T,\ \ n\geq1.$$
(i) and (iv) can be proved using Lemma 4.1 and Definition 4.3, immediately. (vi) can be proved using (vi) in Lemma 4.1, (iii) in Lemma 4.2 and Definition 4.3, immediately. By (4.3), we can get
$$|{\cal{E}}_{\tau,\tau}[X]
-X|\leq|{\cal{E}}_{\tau,\tau}^{\mu,\phi}[X]-X|+|{\cal{E}}^{\mu,\phi}_{\tau,\tau}[X]-X|=0,\ \ P-a.s.$$
Then (ii) is true. Now, we prove (iii). For $\delta\in{\cal{T}}^0_{0,T},$ let $\{\rho_n\}_{n\geq1}\subset{\cal{T}}^0_{0,T}$ be a decreasing sequence such that $\rho_n\leq \delta$ and $\|\rho_n-\rho\|_\infty\rightarrow 0,$ as $n\rightarrow\infty.$ By (iii) in Lemma 4.1, for $ X\in L^\infty({\mathcal {F}}_\delta),$ we have
$${\cal{E}}_{\sigma,\rho_n}[{\cal{E}}_{\rho_n,\delta}[X]]={\cal{E}}_{\sigma,\delta}[X],\ \ P-a.s.\eqno(4.4)$$
By (vi) in Lemma 4.1 and Lemma 2.6, we have ${\cal{E}}_{\delta\wedge\cdot,\delta}[X]\in{\cal{D}}^\infty_{{\cal{F}}}(0,\delta).$ By this and dominated convergence theorem, we have
$$\lim_{n\rightarrow\infty}E\left[|{\cal{E}}_{\rho_n,\delta}[X]-{\cal{E}}_{\rho,\delta}[X]|^2\right]=0.\eqno(4.5)$$
Since
\begin{eqnarray*} &&|{\cal{E}}_{\sigma,\rho_n}[{\cal{E}}_{\rho_n,\delta}[X]]-{\cal{E}}_{\sigma,\rho}[{\cal{E}}_{\rho,\delta}[X]]|\\
&\leq&|{\cal{E}}_{\sigma,\rho_n}[{\cal{E}}_{\rho_n,\delta}[X]]-{\cal{E}}_{\sigma,\rho_n}[{\cal{E}}_{\rho,\delta}[X]]|
+|{\cal{E}}_{\sigma,\rho_n}[{\cal{E}}_{\rho,\delta}[X]]-{\cal{E}}_{\sigma,\rho}[{\cal{E}}_{\rho,\delta}[X]]|,\ \ P-a.s.
\end{eqnarray*}
Thus by (4.5), (ii) in Lemma 4.2 and Definition 4.3, we can get
$$\lim_{n\rightarrow\infty}E[|{\cal{E}}_{\sigma,\rho_n}[{\cal{E}}_{\rho_n,\delta}[X]]
-{\cal{E}}_{\sigma,\rho}[{\cal{E}}_{\rho,\delta}[X]]|^2]=0.$$
By this and (4.4), we have ${\cal{E}}_{\sigma,\rho}[{\cal{E}}_{\rho,\delta}[X]]={\cal{E}}_{\sigma,\delta}[X].$ Thus, we have
$${\cal{E}}_{\sigma,\rho}[{\cal{E}}_{\rho,\tau_n}[X]]={\cal{E}}_{\sigma,\tau_n}[X],\ \ P-a.s.\eqno(4.6)$$
By Definition 4.3, we have
$$\lim_{n\rightarrow\infty}\|{\cal{E}}_{\rho,\tau_n}[X]-{\cal{E}}_{\rho,\tau}[X]\|_\infty=0.$$
From this, (4.3) and the same proof of (i) in Lemma 4.2, we can get
$$\lim_{n\rightarrow\infty}\|{\cal{E}}_{\sigma,\rho}[{\cal{E}}_{\rho,\tau_n}[X]]
-{\cal{E}}_{\sigma,\rho}[{\cal{E}}_{\rho,\tau}[X]]\|_\infty= 0.\eqno(4.7)$$
By (4.6), (4.7) and Definition 4.3, we have
$${\cal{E}}_{\sigma,\rho}[{\cal{E}}_{\rho,\tau}[X]]={\cal{E}}_{\sigma,\tau}[X],\ \ P-a.s.$$
Thus, (iii) is true. By (v) in Lemma 4.1 and the similar argument as (i)-(iv), we can obtain (v).
The proof is complete.\ \  $\Box$\\

\textbf{Step III.}   For $\tau\in{\cal{T}}_{0,T},$ we denote the following space: $\widehat{{\cal{D}}}^2_{{\cal{F}}}(0,\tau)=\{K\in {\cal{D}}^2_{{\cal{F}}}(0,\tau);$ there exists $K^n_{\tau\wedge t}=\int_0^{\tau\wedge t}\gamma_s^nds$ with $\gamma_s^n\in L^\infty_{{\cal{F}}}(0,\tau), n\geq1,$ such that, $K^n\rightarrow K$ in $L^2_{\cal{F}}(0,\tau)$ and for each $t\in[0,T],$ $K^n_{\tau\wedge t} \rightarrow K_{\tau\wedge t}$ in $L^2({\cal{F}}_T),$ as $n\rightarrow\infty\}.$

Now, let $\tau\in{\cal{T}}_{0,T}$ $X\in L^2 ({\cal{F}}_{\tau})$ and $X_n=(X\vee(-n))\wedge n, n\geq1.$ Clearly, $X_n\in L^\infty ({\cal{F}}_{\tau})$ and $X^n\rightarrow X$ in $L^2({\cal{F}}_T),$ as $n\rightarrow\infty.$ For $K\in \widehat{{\cal{D}}}^2_{{\cal{F}}}(0,\tau)$, let $K^n_{\tau\wedge t}=\int_0^{\tau\wedge t}\gamma_s^nds$ with $\gamma_s^n\in L^\infty_{{\cal{F}}}(0,\tau), n\geq1,$ such that, $K^n\rightarrow K$ in $L^2_{\cal{F}}(0,\tau)$ and for each $t\in[0,T],$ $K^n_{\tau\wedge t} \rightarrow K_{\tau\wedge t}$ in $L^2({\cal{F}}_T),$ as $n\rightarrow\infty.$ Consequently, by (4.3), we have
\begin{eqnarray*}
\ \ \ \  &&E|{\cal{E}}_{\tau\wedge t,\tau}[X_n;K^n]-{\cal{E}}_{\tau\wedge t,\tau}[X_m;K^m]|^2\\
&\leq&2E|{\cal{E}}_{\tau\wedge t,\tau}^{\mu,\phi}[X_n-X_m;K^n-K^m]|^2+2E|{\cal{E}}_{\tau\wedge t,\tau}^{-\mu,-\phi}[X_n-X_m;K^n-K^m]|^2\\
&\leq&4E|{\cal{E}}_{\tau\wedge t,\tau}^{\mu,\phi}[X_n-X_m;K^n-K^m]+K^n_{\tau\wedge t}-K^m_{\tau\wedge t}|^2
\\&&+4E|{\cal{E}}_{\tau\wedge t,\tau}^{-\mu,-\phi}[X_n-X_m;K^n-K^m]+K^n_{\tau\wedge t}-K^m_{\tau\wedge t}|^2+8E|K^n_{\tau\wedge t}-K^m_{\tau\wedge t}|^2.
\end{eqnarray*}
By this and Lemma 2.5, we have for $t\in[0,T],$ $\{{\cal{E}}_{\tau\wedge t,\tau}[X^n;K^n]\}_{n\geq1}$ is a Cauchy sequence in $L^2 ({\cal{F}}_T).$ For $t\in[0,T],$ we define
$${\cal{E}}_{\tau\wedge t,\tau}[X;K]=\lim_{n\rightarrow\infty}{\cal{E}}_{\tau\wedge t,\tau}[X^n;K^n]\ \ \textmd{in}\ \ L^2 ({\cal{F}}_T). \eqno(4.8)$$
By (4.3), (4.8) and Lemma 2.5, for $X, X'\in L^2({\mathcal {F}}_\tau)$ and $K, K'\in\widehat{{\cal{D}}}^2_{{\cal{F}}}(0,\tau),$ we can get $\forall t\in[0,T],$
$${\cal{E}}^{-\mu,-\phi}_{{\tau\wedge t},\tau}[X-X';K-K']\leq{\cal{E}}_{{\tau\wedge t},\tau}[X;K]-{\cal{E}}_{{\tau\wedge t},\tau}[X';K']
\leq{\cal{E}}_{{\tau\wedge t},\tau}^{\mu,\phi}[X-X';K-K'],\ \ P-a.s.$$
From this and the same proof of Proposition 3.6, it follows that
\begin{eqnarray*}
\lim_{n\rightarrow\infty}E\left[\sup_{t\in[0,T]}|{\cal{E}}_{\tau\wedge t,\tau}[X_n;K^n]+K^n_{\tau\wedge t}-{\cal{E}}_{\tau\wedge t,\tau}[X;K]-K_{\tau\wedge t}|\right]^2=0.
\end{eqnarray*}
From this, ${\cal{E}}_{\tau\wedge t,\tau}[X;K]+K_{\tau\wedge t}$ is RCLL. Thus ${\cal{E}}_{\tau\wedge t,\tau}[X;K]$ is RCLL. By this and (4.8), we can give the following Definition 4.5.\\\\
\textbf{Definition 4.5 } Let ${\cal{F}}$-evaluation ${\cal{E}}_{s,t}[\cdot]$ satisfy (H1) and (H2), $\sigma, \tau\in{\cal{T}}_{0,T}$ and $\sigma\leq\tau,$ $K\in \widehat{{\cal{D}}}^2_{{\cal{F}}}(0,\tau)$ and $X\in L^2 ({\cal{F}}_{\tau}).$ For each $t\in[0,T],$ we set $\eta_{\tau\wedge t}:={\cal{E}}_{\tau\wedge t,\tau}[X;K].$ Then we define
$${\cal{E}}_{\sigma,\tau}[X;K]:=\eta_\sigma \ \ \ \textmd{and}\ \ \ {\cal{E}}_{\sigma,\tau}[X]:={\cal{E}}_{\sigma,\tau}[X;0]. $$

Now, we have extended the definition of ${\cal{F}}$-evaluation ${\cal{E}}_{s,t}[\cdot;K]$ to with $\sigma, \tau\in{\cal{T}}_{0,T}$ for squared integrable terminal variable and a very special $K.$ Moreover, we have\\\\
\textbf{Lemma 4.6} \textit{The system of operators
$${\cal{E}}_{\sigma,\tau}[\cdot]:\ L^2({\cal{F}}_\tau)\longrightarrow L^2({\cal{F}}_\sigma),\ \ \sigma\leq\tau,\ \sigma\in{\cal{T}}_{0,T},\ \tau\in{\cal{T}}_{0,T},$$
satisfy}

\textit{(i) Monotonicity: ${\cal{E}}_{\sigma,\tau}[\xi]\geq{\cal{E}}_{\sigma,\tau}[\eta], P-a.s., $ if $\xi, \eta \in L^2({\cal{F}}_\tau)$ and $\xi\geq \eta,\ P-a.s.;$}

\textit{(ii)  ${\cal{E}}_{\tau,\tau}[\xi]=\xi, P-a.s., $ if $\xi\in L^2({\mathcal {F}}_\tau);$}

\textit{(iii) Consistency: ${\cal{E}}_{\sigma,\rho}[{\cal{E}}_{\rho,\tau}[\xi]]={\cal{E}}_{\sigma,\tau}[\xi], P-a.s.,$ if $\sigma\leq \rho\leq \tau$ and $\xi\in L^2({\mathcal {F}}_\tau), \ \rho\in{\cal{T}}_{0,T};$}

\textit{(iv) "0-1 Law":  $1_A{\cal{E}}_{\sigma,\tau}[\xi]={\cal{E}}_{\sigma,\tau}[1_A\xi], P-a.s., $ if $A\in{\mathcal {F}}_\sigma,\ \xi\in L^2({\mathcal {F}}_\tau);$}

\textit{(v) For $\tau\in{\cal{T}}_{0,T},$ $K\in \widehat{{\cal{D}}}^2_{{\cal{F}}}(0,\tau),$ $${\cal{E}}_{\sigma,\tau'}[\cdot;K]:\ L^2({\cal{F}}_{\tau'})\longrightarrow L^2({\cal{F}}_\sigma),\ \ \sigma\leq\tau'\leq\tau,\ \sigma,\tau'\in{\cal{T}}_{0,T},$$ satisfies the above (i)-(iii) and}
$$1_A{\cal{E}}_{\sigma,\tau'}[\xi;K]=1_A{\cal{E}}_{\sigma,\tau'}[1_A\xi;K],\ \ P-a.s. \ \ \textrm{if} \ \ A\in{\mathcal {F}}_\sigma,\ \xi\in L^2({\mathcal {F}}_{\tau'});$$

\textit{(vi) For $\tau\in{\cal{T}}_{0,T},$ $K\in \widehat{{\cal{D}}}^2_{{\cal{F}}}(0,\tau)$ and $\xi\in L^2({\mathcal {F}}_\tau),$ ${\cal{E}}_{\tau\wedge \cdot,\tau}[\xi;K]$ is RCLL and for $X, X'\in L^2({\mathcal {F}}_\tau)$ and $K, K'\in\widehat{{\cal{D}}}^2_{{\cal{F}}}(0,\tau),$ we have}
$${\cal{E}}^{-\mu,-\phi}_{\sigma,\tau}[X-X';K-K']\leq{\cal{E}}_{\sigma,\tau}[X;K]-{\cal{E}}_{\sigma,\tau}[X';K']
\leq{\cal{E}}_{\sigma,\tau}^{\mu,\phi}[X-X';K-K'],\ \ P-a.s.$$
\emph{Proof.} Clearly, we only need prove (v) and (vi). Given $\tau\in{\cal{T}}_{0,T},$ for $\sigma,\tau'\in{\cal{T}}_{0,T}$ and $\sigma\leq\tau'\leq\tau,$ we can firstly prove (vi) and that ${\cal{E}}_{\sigma,\tau'}[\cdot;K]$ satisfies (i) by Lemma 4.4 and Definition 4.5, immediately. Then we can prove that ${\cal{E}}_{\sigma,\tau'}[\cdot;K]$ satisfies (ii) by (vi) like the proof of (ii) in Lemma 4.4. In the following, we will prove ${\cal{E}}_{\sigma,\tau'}[\cdot;K]$ satisfies (iii). For $K\in \widehat{{\cal{D}}}^2_{{\cal{F}}}(0,\tau),$  let $K^n_{\tau\wedge t}=\int_0^{\tau\wedge t}\gamma_s^nds$ with $\gamma_s^n\in L^\infty_{{\cal{F}}}(0,\tau),$ such that, $K^n\rightarrow K$ in $L^2_{\cal{F}}(0,\tau)$ and for each $t\in[0,T],$ $K^n_{\tau\wedge t} \rightarrow K_{\tau\wedge t}$ in $L^2({\cal{F}}_T),$ as $n\rightarrow\infty.$ For $X\in L^2 ({\cal{F}}_{\tau'}),$ let $X_n=(X\vee(-n))\wedge n.$ For $\rho\in{\cal{T}}_{0,T}$ and $\sigma\leq\rho\leq\tau'.$
by (vi), comparison theorem and "Consistency", we have $P-a.s.,$
\begin{eqnarray*}
&&{\cal{E}}_{\sigma,\rho}[{\cal{E}}_{\rho,\tau'}[X^n;K^n];K^n]+K^n_\sigma
-{\cal{E}}_{\sigma,\rho}[{\cal{E}}_{\rho,\tau'}[X;K];K]-K_\sigma\\
&\leq&{\cal{E}}^{\mu,\phi}_{\sigma,\rho}[{\cal{E}}_{\rho,\tau'}[X^n;K^n]
-{\cal{E}}_{\rho,\tau'}[X^n;K^n];K^n-K]+K^n_\sigma-K_\sigma\\
&\leq&{\cal{E}}^{\mu,\phi}_{\sigma,\rho}[{\cal{E}}^{\mu,\phi}_{\rho,\tau'}[X^n-X;K^n-K];K^n-K]+K^n_\sigma-K_\sigma\\
&=&{\cal{E}}^{\mu,\phi}_{\sigma,\tau'}[X^n-X;K^n-K]+K^n_\sigma-K_\sigma.
\end{eqnarray*}
Similarly, we have $P-a.s.,$
\begin{eqnarray*}
&&{\cal{E}}_{\sigma,\rho}[{\cal{E}}_{\rho,\tau'}[X^n;K^n];K^n]+K^n_\sigma
-{\cal{E}}_{\sigma,\rho}[{\cal{E}}_{\rho,\tau'}[X;K];K]-K_\sigma\\
&\geq&{\cal{E}}^{-\mu,-\phi}_{\sigma,\tau'}[X^n-X;K^n-K]+K^n_\sigma-K_\sigma.
\end{eqnarray*}
Thus, by the above two inequalities and Lemma 2.5, we have
$${\cal{E}}_{\sigma,\rho}[{\cal{E}}_{\rho,\tau'}[X^n;K^n];K^n]+K^n_\sigma
\rightarrow{\cal{E}}_{\sigma,\rho}[{\cal{E}}_{\rho,\tau'}[X;K];K]+K_\sigma,\ \ \textmd{in}\ \  L^2 ({\cal{F}}_T),$$
as $n\rightarrow\infty.$ Similar argument as the above gives
$${\cal{E}}_{\sigma,\tau'}[X^n;K^n]+K^n_\sigma
\rightarrow{\cal{E}}_{\sigma,\tau'}[X;K]+K_\sigma,\ \ \textmd{in}\ \  L^2 ({\cal{F}}_T),\eqno(4.9)$$
as $n\rightarrow\infty.$ By (v) in Lemma 4.4, we have
$${\cal{E}}_{\sigma,\rho}[{\cal{E}}_{\rho,\tau'}[X^n;K^n];K^n]={\cal{E}}_{\sigma,\tau'}[X^n;K^n],\ \ P-a.s.$$
From the above three equalities, it follows that $${\cal{E}}_{\sigma,\rho}[{\cal{E}}_{\rho,\tau'}[X;K];K]={\cal{E}}_{\sigma,\tau'}[X;K],\ \ P-a.s.$$
Thus ${\cal{E}}_{\sigma,\tau'}[\cdot;K]$ satisfies (iii). By (4.9), for $A\in{\mathcal {F}}_\sigma,$  we have
$$1_A{\cal{E}}_{\sigma,\tau'}[X^n;K^n]+1_AK^n_\sigma\rightarrow1_A{\cal{E}}_{\sigma,\tau'}[X;K]+1_AK_\sigma,\ \ \textmd{in}\ \  L^2 ({\cal{F}}_T),$$
and
$$1_A{\cal{E}}_{\sigma,\tau'}[1_AX^n;K^n]+1_AK^n_\sigma\rightarrow1_A{\cal{E}}_{\sigma,\tau'}[1_AX;K]+1_AK_\sigma,\ \ \textmd{in}\ \  L^2 ({\cal{F}}_T),$$
as $n\rightarrow\infty.$ Thus, by (v) in Lemma 4.4, we have
$$1_A{\cal{E}}_{\sigma,\tau'}[X;K]=1_A{\cal{E}}_{\sigma,\tau'}[1_AX;K],\ \ P-a.s.$$
The proof is complete.\ \ $\Box$\\

The following Lemma 4.7 is an optional stopping theorem for locally bounded ${\cal{E}}_{s,t}[\cdot;K]$-supermartingales, which is crucial in the proof of Lemma 4.8 and Proposition 5.5.\\\\
\
\textbf{Lemma 4.7}  \textit{Let ${\cal{F}}$-evaluation ${\cal{E}}_{s,t}[\cdot]$ satisfy (H1) and (H2), $K_t=\int_0^t\gamma_sds$ with $\gamma_s\in L^\infty_{{\cal{F}}}(0,T),$ $\tau\in{\cal{T}}_{0,T}$ and $Y\in {\cal{D}}^2_{{\cal{F}}}(0,T)$ be an ${\cal{E}}_{s,t}[\cdot;K]$-supermartingale (resp. ${\cal{E}}_{s,t}[\cdot;K]$-submartingale) with $Y\in {\cal{D}}^\infty_{{\cal{F}}}(0,\tau)$ and $Y_\tau\in L^\infty({\cal{F}}_\tau).$ Then for $\sigma,\ \tau'\in{\cal{T}}_{0,T}$ satisfing $\sigma\leq\tau'\leq\tau,$ we have}
$${\cal{E}}_{\sigma,\tau'}[Y_{\tau'};K]\leq Y_\sigma\ (resp. \geq),\ \ P-a.s.$$
\emph{Proof.}  We only prove the ${\cal{E}}_{s,t}[\cdot;K]$-supermartingales case. The ${\cal{E}}_{s,t}[\cdot;K]$-submartingales case is similar. we prove it by two steps.

Step A. Let $\sigma\in{\cal{T}}_{0,T},\tau'\in{\cal{T}}^0_{0,T},\sigma\leq\tau',$ $K'\in {\cal{D}}^2_{{\cal{F}}}(0,T)$ and $Y'\in {\cal{D}}^2_{{\cal{F}}}(0,T)$ be an ${\cal{E}}_{s,t}[\cdot;K']$-supermartingale. Let $\{\sigma_n\}_{n\geq1}\subset{\cal{T}}^0_{0,T}$ satisfy $\sigma_n\leq\tau'$ and $\sigma_n\searrow\sigma,$ as $n\rightarrow\infty.$ By Lemma 4.1, we can get ${\cal{E}}_{\sigma_n,\tau'}[\cdot;K']$ satisfy (i)-(iii) in Lemma 4.1 and (4.1). Thus by the proof of Peng [14, Lemma 10.10], we can get ${\cal{E}}_{\sigma_n,{\tau'}}[Y'_{\tau'};K']\leq Y'_{\sigma_n}.$ By the right continuity of ${\cal{E}}_{\tau'\wedge t,\tau'}[Y'_{\tau'};K']$ and $Y'$, we have ${\cal{E}}_{{\sigma},{\tau'}}[Y'_{\tau'};K']\leq Y'_{\sigma},\ \ P-a.s.$

Step B.  Let $\sigma, \tau'\in{\cal{T}}_{0,T}, \sigma\leq\tau'\leq\tau,$
and $\{\tau'_n\}_{n\geq1}\subset{\cal{T}}^0_{0,T}$ is a decreasing sequence such that $\|\tau'_n-\tau'\|_\infty\rightarrow 0.$ By Step A, we have
$${\cal{E}}_{\sigma,\tau'_n}[Y_{\tau'_n};K]\leq Y_\sigma,\ \ P-a.s.\eqno(4.10)$$
Since
\begin{eqnarray*}
\ \ \ \ \ \ \ \ \ \ \ \ \ \ &&|{\cal{E}}_{\sigma,\tau'_n}[Y_{\tau'_n};K]-{\cal{E}}_{\sigma,\tau'}[Y_{\tau'};K]|\\
&\leq&|{\cal{E}}_{\sigma,\tau'_n}[Y_{\tau'_n};K]-{\cal{E}}_{\sigma,\tau'_n}[Y_{\tau'};K]|
+|{\cal{E}}_{\sigma,\tau'_n}[Y_{\tau'};K]-{\cal{E}}_{\sigma,\tau'}[Y_{\tau'};K]|,\ \ \ \ \ \ \ \ \ \ \ \ \ \ \ \ \ (4.11)
\end{eqnarray*}
and $Y_{\tau'_n}\rightarrow Y_{\tau'}$ in $L^2({\cal{F}}_T)$ as $n\rightarrow\infty,$ thus by (ii) in Lemma 4.2 and Definition 4.3, we have
$$\lim_{n\rightarrow\infty}E[|{\cal{E}}_{\sigma,\tau'_n}[Y_{\tau'_n};K]-{\cal{E}}_{\sigma,\tau'}[Y_{\tau'};K]|^2]
=0. \eqno(4.12)$$
By (4.10) and (4.12), we complete this proof. \ \ $\Box$\\\\
\
\textbf{Lemma 4.8} \textit{Let $g$ satisfy (A1) and (A2), $K_t=\int_0^t\gamma_sds$ with $\gamma_s\in L^\infty_{{\cal{F}}}(0,T),$ $\tau\in{\cal{T}}_{0,T}$ and $Y\in {\cal{D}}^2_{{\cal{F}}}(0,T)$ be an ${\cal{E}}^g_{s,t}[\cdot;K]$-supermartingale with $Y\in {\cal{D}}^\infty_{{\cal{F}}}(0,\tau)$ and $Y_\tau\in L^\infty({\cal{F}}_\tau).$  Then there exists a process $A_s\in{\mathcal{D}}^2_{\cal{F}}(0,\tau)$, which is increasing with $A_0=0,$ such that for $ \sigma,\ \tau'\in{\cal{T}}_{0,T}$ satisfying $\sigma\leq\tau'\leq\tau,$ we have}
$$Y_\sigma={\cal{E}}^g_{\sigma,\tau'}[Y_{\tau'};K+A],\ \ P-a.s.$$
\textit{Proof.} By Remark 3.3 and the above arguments of this section, we can get the optimal stopping theorem (Lemma 4.7) also holds true for $Y_t$. That is, for $\sigma,\ \tau'\in{\cal{T}}_{0,T}$ satisfying $\sigma\leq\tau'\leq\tau,$ we have
$${\cal{E}}^g_{\sigma,\tau'}[Y_{\tau'};K]\leq Y_\sigma,\ \ P-a.s.\eqno(4.13)$$
By (i) in Remark 2.2 and (4.13), for $\sigma,\ \tau'\in{\cal{T}}_{0,T}$ satisfying $\sigma\leq\tau'\leq\tau,$ we have
$${\cal{E}}^{g^K}_{\sigma,{\tau'}}[Y_{\tau'}+K_{\tau'}]={\cal{E}}^g_{\sigma,{\tau'}}[Y_{\tau'};K]+K_\sigma\leq Y_\sigma+K_\sigma,\ \ P-a.s.$$
By this, we can obtain a result similar as Peng [16, Lemma 3.8] by a similar argument. Then by the similar proof as Peng [15, Theorem 3.3] or Peng [16, Theorem 3.9], we can get that there exists $A\in {\cal{D}}^2_{{\cal{F}}}(0,\tau)$ such that for $\sigma,\ \tau'\in{\cal{T}}_{0,T}$ satisfying $\sigma\leq\tau'\leq\tau,$  we have
$$Y_\sigma+K_\sigma={\cal{E}}^{g^K}_{\sigma,{\tau'}}[Y_\tau+K_{\tau'};A],\ \ P-a.s.$$
From this, we can get $Y_\sigma={\cal{E}}^g_{\sigma,{\tau'}}[Y_{\tau'};K+A],\ P-a.s.$ The proof is complete. \ \ $\Box$\\

Now, we give the following Lemma 4.9, which is important in the proof of Theorem 5.4. \\\\
\textbf{Lemma 4.9}\textit{ Let ${\cal{F}}$-expectation ${\cal{E}}_{s,t}[\cdot]$ satisfy (H1) and (H2), $K_t=\int_0^t\gamma_sds$ with $\gamma_s\in L^\infty_{{\cal{F}}}(0,T).$ Let $\tau\in{\cal{T}}_{0,T}$ and $X\in L^\infty({\cal{F}}_\tau).$ For $\sigma\in{\cal{T}}_{0,T}$ satisfying $\sigma\leq\tau,$ we set
$$Y_{\sigma}^{\tau,X,K}:={\cal{E}}_{\sigma,\tau}[X;K].$$
Then there exists a pair $(g_s^{\tau,X,K}, Z_s^{\tau,X,K})$ in $L^2_{{\cal{F}}}(0,\tau)\times L^2_{{\cal{F}}}(0,\tau;{\textbf{R}}^d)$ such that $\forall t\in[0,\tau],$
$$|g_{t}^{\tau,X,K}|\leq \mu|Y_{t}^{\tau,X,K}|+\phi(|Z_{t}^{\tau,X,K}|),\ \ P-a.s.$$ and $\forall t\in[0,T],$
$$Y_{\tau\wedge t}^{\tau,X,K}=X+K_\tau-K_{\tau\wedge t}+\int_{\tau\wedge t}^\tau g_r^{\tau,X,K}dr-\int_{\tau\wedge t}^\tau Z_r^{\tau,X,K}dB_r,\ \ P-a.s.$$
Moreover, for $\tau'\in{\cal{T}}_{0,T},$ $X'\in L^\infty({\cal{F}}_{\tau'})$ and $K'_t=\int_0^t\gamma'_sds$ with $\gamma'_s\in L^\infty_{{\cal{F}}}(0,T),$ we have $ \forall t\in[0,\tau\wedge \tau'],$}
$$|g_{t}^{\tau,X,K}-g_{t}^{\tau',X',K'}|\leq  \mu(|Y_{t}^{\tau,X,K}-Y_{ t}^{\tau',X',K'}|)+\phi(|Z_{t}^{\tau,X,K}-Z_{t}^{\tau',X',K'}|),\ \ P-a.s.$$
\emph{Proof.} By (vi) in Lemma 4.4 and "Consistency", for $\sigma,\tau'\in{\cal{T}}_{0,T}$ satisfying $\sigma\leq\tau'\leq\tau,$ we have
$${\cal{E}}^{-\mu,-\phi}_{\sigma,\tau'}[Y_{\tau'}^{\tau,X,K};K]\leq{\cal{E}}_{\sigma,{\tau'}}[Y_{\tau'}^{\tau,X,K};K]
={\cal{E}}_{\sigma,{\tau'}}[{\cal{E}}_{\tau',\tau}[X;K];K]={\cal{E}}_{\sigma,\tau}[X;K]=Y_{\sigma}^{\tau,X,K}.\eqno(4.14)$$
Clearly, one can find the proof of Lemma 4.8 is based on (4.13). Thus, by (4.14), we can get there exists a process $A_s^-\in{\mathcal{D}}^2_{\cal{F}}(0,\tau)$, which is increasing with $A_0^-=0,$ such that for each $t\in[0,T],$ we have
$$Y^{\tau,X,K}_{\tau\wedge t}={\cal{E}}^{-\mu,-\phi}_{{\tau\wedge t},\tau}[X;K+A^-],\ \ P-a.s.\eqno(4.15)$$
Similarly, we also can show there exists a process $A_s^+\in{\mathcal{D}}^2_{\cal{F}}(0,\tau)$, which is increasing with $A_0^+=0$ such that  for each $t\in[0,T],$ we have
$$Y^{\tau,X,K}_{\tau\wedge t}={\cal{E}}^{\mu,\phi}_{{\tau\wedge t},\tau}[X;K-A^+],\ \ P-a.s.\eqno(4.16)$$
By (4.15) and (4.16), we can complete the proof by the similar argument of Peng [14, Proposition 6.6 and Corollary 6.7]. We omit it here. \ \ $\Box$ \\\\
\textbf{Remark 4.10} Let ${\cal{F}}$-expectation ${\cal{E}}_{s,t}[\cdot]$ satisfy (H1) and (H2), $K_t=\int_0^t\gamma_sds$ with $\gamma_s\in L^\infty_{{\cal{F}}}(0,T),$ $\tau\in{\cal{T}}_{0,T}.$ Then for $X\in L^\infty({\cal{F}}_\tau),$ we can get ${\cal{E}}_{\tau\wedge\cdot,\tau}[X;K]\in{\cal{S}}^\infty_{{\cal{F}}}(0,\tau),$ from (4.3), Lemma 2.6 and Lemma 4.9.

%%%%%%%%%%%%%%%%%%%%%%%%%%%%%%%%%%%%%%%%%%%%%%%%%%%%%%%%%%%%%%%%
\section{Doob-Meyer decomposition of ${\cal{E}}_{s,t}[\cdot]$-supermartingales}
%%%%%%%%%%%%%%%%%%%%%%%%%%%%%%%%%%%%%%%%%%%%%%%%%%%%%%%%%%%%%%%%
In this section, we will study the Doob-Meyer decomposition of ${\cal{E}}_{s,t}[\cdot]$-supermartingales. It is obtained in a locally bounded case. Given a function $f: \Omega \times [0,T]\times \mathbf{R}\longmapsto \mathbf{R},$ in this paper, we always suppose $f$ satisfy the following Lipschitz condition:
$$\exists \lambda\geq0,\ s.t.\ |f(t,y_1)-f(t,y_2)|\leq \lambda|y_1-y_2|, \ \forall y_1,\ y_2\in{\mathbf{R}},\ \forall t\in[0,T].$$
Now, we consider the following BSDE denoted by ${\cal{E}}(f,X,T)$ under ${\cal{F}}$-evaluation ${\cal{E}}_{s,t}[\cdot]:$
$$y_s={\cal{E}}_{s,T}\left[X; \int_0^\cdot f(r,y_r)dr\right],\ \ s\in[0,T].$$
\textbf{Theorem 5.1} \textit{Let ${\cal{F}}$-evaluation ${\cal{E}}_{s,t}[\cdot]$ satisfy (H1)  and (H2), $X\in L^\infty({\cal{F}}_T)$ and $f(\cdot,0)\in L^\infty_{{\cal{F}}}(0,T).$ Then ${\cal{E}}(f,X,T)$ has a unique solution $y_t\in {\cal{S}}_{\cal{F}}^\infty (0,T).$ }\\\\
\emph{Proof.} For $y_s\in {\cal{S}}^\infty_{{\cal{F}}}(0,T),$ set $$I(y_s):={\cal{E}}_{s,T}\left[X;\int_0^\cdot f(r,y_r)dr\right],$$
Since $f$ satisfies Lipschitz condition, $y_s\in {\cal{S}}^\infty_{{\cal{F}}}(0,T)$ and $f(\cdot,0)\in L^\infty_{{\cal{F}}}(0,T)$, thus we have
\begin{eqnarray*}
\|f(r,y_r)\|_{L^\infty_{\cal{F}}(0,T)}\leq\|f(r,0)\|_{L^\infty_{\cal{F}}(0,T)}+\lambda\|y_r\|_{L^\infty_{\cal{F}}(0,T)}
<\infty.
\end{eqnarray*}
Then by Remark 4.10, we have $I(y_s)\in {\cal{S}}_{\cal{F}}^\infty (0,T).$ Thus
$$I(\cdot):{\cal{S}}_{\cal{F}}^\infty (0,T)\longmapsto {\cal{S}}_{\cal{F}}^\infty (0,T).$$
By (4.3), for each $y_s^1,\ y_s^2\in {\cal{S}}_{\cal{F}}^\infty (0,T),$ we have
\begin{eqnarray*}
  &&|I(y_s^1)-I(y_s^2)|\\&=&\left|{\cal{E}}_{s,T}\left[X;\int_0^\cdot f(r,y_r^1)dr\right]
  -{\cal{E}}_{s,T}\left[X;\int_0^\cdot f(r,y_r^2)dr\right]\right|\\
  &\leq&\left|{\cal{E}}^{\mu,\phi}_{s,T}\left[0;\int_0^\cdot (f(r,y_r^1)-f(r,y_r^2))dr\right]\right|+\left|{\cal{E}}^{-\mu,-\phi}_{s,T}\left[0;\int_0^\cdot (f(r,y_r^1)- f(r,y_r^2))dr\right]\right|.
\end{eqnarray*}
By Lemma 2.6, we can get
\begin{eqnarray*}
  \left\|{\cal{E}}^{\mu,\phi}_{s,T}\left[0;\int_0^\cdot (f(r,y_r^1)-f(r,y_r^2))dr\right]\right\|_{L^\infty_{\cal{F}}(0,T)}
  &\leq&Te^{\mu T}\left\|f(s,y_s^1)-f(s,y_s^2)\right\|_{L_{\cal{F}}^\infty(0,T)}\\
  &\leq&\lambda T e^{\mu T}\left\|y_s^1-y_s^2\right\|_{L^\infty_{\cal{F}}(0,T)}.
\end{eqnarray*}
Similarly, we have
\begin{eqnarray*}
  \left\|{\cal{E}}^{-\mu,-\phi}_{s,T}\left[0;\int_0^\cdot (f(r,y_r^1)-f(r,y_r^2))dr\right]\right\|_{L^\infty_{\cal{F}}(0,T)}
  \leq\lambda T e^{\mu T}\left\|y_s^1-y_s^2\right\|_{L^\infty_{\cal{F}}(0,T)}.
\end{eqnarray*}
Thus from above three inequalities, there exists a constant $\beta>0$ such that if $T\leq\beta,$ we have
$$\left\|I(y_s^1)-I(y_s^2)\right\|_{L^\infty_{\cal{F}}(0,T)}\leq\frac{1}{2}\left\|y_s^1-y_s^2\right\|_{L^\infty_{\cal{F}}(0,T)}.$$
Consequently, in the case that $T\leq\beta,$ $I(\cdot)$ is a strict contraction. The proof is complete.

In the case that $T>\beta$, we can complete the proof using a "patching-up" method given in Hu et al. [7, Proposition 4.4]. We take a partition of $[0,T]:\ 0=t_0<t_1<\cdots<t_N=T$ such that $\max_n|t_n-t_{n-1}|\leq\beta$. In view of Lemma 2.6, we can prove ${\cal{E}}(f,t_N,X)$ has a unique solution on $[t_{N-1},t_N]$ by the above argument, we denote the solution by $y^N_s,$ $s\in[t_{N-1},t_N].$ Similarly, we can solve ${\cal{E}}(f,t_{n-1},y_{t_{n-1}}^n)$ on $[t_{n-2},t_{n-1}],$ and denote its solution by $y^{n-1}_s,$ $s\in[t_{n-2},t_{n-1}],$ $2\leq n\leq N$. Now, we set $y_s:=y^{n}_s,\ s\in [t_{n-1},t_{n}],\ 1\leq n\leq N,$ we will show $y_t$ is a solution of ${\cal{E}}(f,T,X)$ on $[0,T].$

Clearly,  $y_s$ is a solution of ${\cal{E}}(f,T,X)$ on $[t_{N-1},T].$ Assuming $y_s$ is a solution of ${\cal{E}}(f,T,X)$ on $[t_{m},T],\ 1<m\leq N-1,$  then by above settings and  "Consistency" of ${\cal{E}}$, for $s\in[t_{m-1},t_{m}],$ we have
\begin{eqnarray*}
  y_s= y_s^m
  &=&{\cal{E}}_{s,t_m}\left[y_{t_{m}}^m;\int_0^\cdot f(r,y_r)dr\right]\\
  &=&{\cal{E}}_{s,t_m}\left[y_{t_{m}};\int_{0}^\cdot f(r,y_r)dr\right]\\
  &=&{\cal{E}}_{s,t_m}\left[{\cal{E}}_{{t_{m}},T}\left[X;\int_{0}^\cdot f(r,y_r)dr\right]
  ;\int_{0}^\cdot f(r,y_r)dr\right]\\
  &=&{\cal{E}}_{s,T}\left[X;\int_0^\cdot f(r,y_r)dr\right].
\end{eqnarray*}
Thus $y_t$ is also a solution on $[t_{m-1},T].$ By induction, we can get $y_t$ is a solution on $[0,T].$

If $\hat{y}_t\in {\cal{S}}_{\cal{F}}^\infty (0,T)$ is another solution of ${\cal{E}}(f,T,X)$ on $[0,T].$ Clearly by the above argument, we get  $\hat{y}_s=y_s,\ s\in [t_{N-1},N].$  Similarly, we can also get $\hat{y}_s=y_s,\ s\in [t_{n-1},t_n],\ 1\leq n\leq N-1.$ Thus $\hat{y}_s=y_s,\ s\in [0,T].$ The proof is complete.  \ \  $\Box$\\

By the similar arguments as Peng [14, Proposition 7.3 and Corollary 7.4], we can get the following comparison theorem for ${\cal{E}}(f,T,X).$ We omit its proof here.\\\\
\
\textbf{Theorem 5.2}\textit{ Let ${\cal{F}}$-evaluation ${\cal{E}}_{s,t}[\cdot]$ satisfy (H1) and (H2), $X\in L^\infty({\cal{F}}_T),$ $f(\cdot,0)\in L^\infty_{{\cal{F}}}(0,T).$  Let $y_s$ be the solution of ${\cal{E}}(f,T,X)$ and $\bar{y}_s$ be the solution of the following ${\cal{E}}(f+\eta_s,T,\bar{X})$:
$$\bar{y}_s={\cal{E}}_{s,T}\left[\bar{X};\int_0^\cdot(f(r,\bar{y}_r)+\eta_r)dr\right],\ \  t\in[0,T],$$
where $\bar{X}\in L^\infty({\cal{F}}_T)$ and $\eta_t\in L_{\cal{F}}^\infty (0,T)$ satisfy
$$\bar{X}\geq X,\ \ \ \eta_s\geq0,\ \ dP\times dt-a.e.$$
Then we have $\forall s\in[0,T],$}
$$\bar{y}_s\geq y_s,\ \ P-a.s.$$\\
\
\textbf{Remark 5.3 }
\begin{itemize}
\item[(i)] Let ${\cal{F}}$-evaluation ${\cal{E}}_{s,t}[\cdot]$ satisfy (H1)  and (H2). Clearly, if $y_s$ is the solution of ${\cal{E}}(f,T,X),$ then process $y_s$ is an ${\cal{E}}_{s,t}[\cdot;\int_0^\cdot f(r,y_r)dr]$-martingale on $[0,T].$ Thus we can also get that $y_s$ is the unique solution of ${\cal{E}}(f,t,y_t)$ on $[0,t].$
\item[(ii)] Theorem 5.1 and Theorem 5.2 are for ${\cal{E}}(f,T,X)$ with given deterministic terminal time $T.$ In fact, we can also obtain the same conclusion for ${\cal{E}}(f,\tau,X)$ with $\tau\in{\cal{T}}_{0,T},$ from the same arguments.
 \end{itemize}

The following Theorem 5.4 is a Doob-Meyer type decomposition for locally bounded ${\cal{E}}_{s,t}[\cdot;K]$-supermartingales, which generalizes the corresponding result in Lemma 4.9.\\\\
\
\textbf{Theorem 5.4} \textit{ Let ${\cal{F}}$-evaluation ${\cal{E}}_{s,t}[\cdot]$ satisfy (H1) and (H2), $\tau\in{\cal{T}}_{0,T},$ $Y_s\in {\cal{S}}^2_{{\cal{F}}}(0,T)$ is an ${\cal{E}}_{s,t}[\cdot]$-supermartingale with $Y_s\in {\cal{S}}^\infty_{{\cal{F}}}(0,\tau).$ Then there exists a process $A_s\in {\cal{S}}^2_{{\cal{F}}}(0,\tau)$, which is increasing with $A_0=0$ such that} $\forall t\in[0,T],$
$${\cal{E}}_{t\wedge\tau,\tau}[Y_\tau;A]= Y_{t\wedge\tau},\ \ P-a.s.,$$
\emph{and there exists a pair $(g_s, Z_s)$  in $L^2_{{\cal{F}}}(0,\tau)\times L^2_{{\cal{F}}}(0,\tau;{\textbf{R}}^d)$ such that for $t\in[0,\tau],$ }
$$|g_{t}|\leq \mu|Y_{t}|+\phi(|Z_{t}|),\ \ dP\times dt-a.e,$$
and $\forall t\in[0,T],$
$$Y_{\tau\wedge t}=Y_\tau+A_\tau-A_{\tau\wedge t}+\int_{\tau\wedge t}^\tau g_rdr-\int_{\tau\wedge t}^\tau Z_rdB_r,\ \ P-a.s.$$
\emph{Moreover for any ${\cal{E}}_{s,t}[\cdot]$-supermartingale $Y'_s\in {\cal{S}}^2_{{\cal{F}}}(0,T)$ with $Y'_s\in S^\infty_{{\cal{F}}}(0,\tau'),$ the corresponding pair $(g'_s, Z'_s)$ in $L^2_{{\cal{F}}}(0,\tau')\times L^2_{{\cal{F}}}(0,\tau';{\textbf{R}}^d)$ satisfies for $ t\in[0,\tau\wedge \tau'],$}
$$|g_{t}-g'_{t}|\leq  \mu(|Y_{t}-Y'_{t}|)+\phi(|Z_{t}-Z'_{t}|),\ \ dP\times dt-a.e.$$
\emph{Proof.}  For $n\geq1,$ we consider the following BSDE under ${\cal{F}}$-evaluation ${\cal{E}}_{s,t}[\cdot]$:
$$y_{t\wedge\tau}^n={\cal{E}}_{t\wedge\tau,\tau}\left[Y_\tau;\int_0^\cdot n(Y_s-y_s^n)ds\right],\ \ \ t\in[0,T].\eqno(5.1)$$
By Theorem 5.1 and Remark 5.3, the above BSDE (5.1) has a unique solution $y_t^n\in S^\infty_{{\cal{F}}}(0,\tau).$ Then we have the following Proposition 5.5.  \\\\
\textbf{Proposition 5.5} \textit{For $n\geq1$ and each $t\in[0,T],$ we have} $$Y_{t\wedge\tau}\geq y_{t\wedge\tau}^{n+1}\geq y_{t\wedge\tau}^n,\ \ \ P-a.s.$$
\emph{Proof.} With the help of the optional stopping theorem (Lemma 4.7), Theorem 5.1, Theorem 5.2 and Remark 5.3, we can obtain this proposition from the argument of Peng [14, Lemma 8.3], immediately. \ \ $\Box$\\

Set $$A_{t\wedge\tau}^n:=\int_0^{t\wedge\tau}n(Y_s-y_s^n)ds,\ \ t\in[0,T],\ \ n\geq1.\eqno(5.2)$$
By Proposition 5.5, $A_{t\wedge\tau}^n\in{\cal{S}}^\infty_{{\cal{F}}}(0,\tau),$ and is increasing with $A_0=0.$
Then by (5.1) and (5.2), we have $\forall t\in[0,T],$
$$y^n_{t\wedge\tau}={\cal{E}}_{t\wedge\tau,\tau}[Y_\tau;A^n],\ \ P-a.s.\eqno(5.3)$$
Thus by Lemma 4.9, there exists a pair $(g^n_s, Z^n_s)$  in $L^2_{{\cal{F}}}(0,\tau)\times L^2_{{\cal{F}}}(0,\tau;{\textbf{R}}^d)$ such that $\forall t\in[0,\tau],$
$$|g_t^n|\leq\mu|y^n_t|+\phi(|Z^n_t|),\ \ P-a.s.,\ \  n\geq1,\eqno(5.4)$$
$$|g_t^n-g_t^m|\leq\mu|y^n_t-y^m_t|+\phi(|Z^n_t-Z^m_t|),\ \ P-a.s.,\ \ m, n\geq1, \eqno(5.5)$$
and $\forall t\in[0,T],$
$$y_{t\wedge\tau}^n=Y_\tau+A_\tau^n-A_{t\wedge\tau}^n+\int_{t\wedge\tau}^\tau g^n_sds-\int_{t\wedge\tau}^\tau Z^n_sdB_s,\ \ P-a.s.,\ \  n\geq1.\eqno(5.6)$$
Moreover for an ${\cal{E}}_{s,t}[\cdot]$-supermartingale $Y'_s\in {\cal{S}}^2_{{\cal{F}}}(0,T)$ with $Y'_s\in S^\infty_{{\cal{F}}}(0,\tau'),$ the corresponding pair $(g'^n_s, Z'^n_s)$  in $L^2_{{\cal{F}}}(0,\tau')\times L^2_{{\cal{F}}}(0,\tau';{\textbf{R}}^d)$ satisfies $\forall t\in[0,\tau\wedge \tau'],$
$$|{g}_{t}^n-{g'}_{t}^{n}|\leq  \mu(|{y}_{t}^{n}-{y'}_{t}^{n}|)+\phi(|{Z}_{t}^{n}-{Z'}_{t}^{n}|),\ \ P-a.s.,\ \ n\geq1.\eqno(5.7)$$
We further have\\\\
\textbf{Proposition 5.6}\textit{ There exists a constant $C$ independent on $n$, such that}
$$(i)\ \ E\int_0^\tau|Z^n_s|^2ds\leq C\ \ \ and\ \ \ (ii) \ \ E|A^n_\tau|^2\leq C.$$
\emph{Proof.} The proof is similar as Zheng and Li [19, Proposition 4.2], we give it here for convenience. In this proof, $C$  is assumed as a constant independent on $n$, its value may change line by line. By Proposition 5.5, we get that $y_{t\wedge\tau}^1\leq y_{t\wedge\tau}^n\leq y_{t\wedge\tau}^{n+1}\leq Y_{t\wedge\tau}.$ Thus, we have $$\|y^n_t\|_{L_{\cal{F}}^\infty(0,\tau)}\leq C,\ \ n\geq1.\eqno(5.8)$$
By (5.6), (5.4), (5.8) and the fact that $\phi$ has a linear growth, we have
\begin{eqnarray*}
E|A_\tau^n|^2&\leq& 3E|y^n_0-y^n_\tau|^2+3TE\int_0^\tau|g^n_s|^2ds+3E\int_0^\tau|Z^n_s|^2ds\\
&\leq&C+ 3TE\int_0^\tau(\mu|y^n_s|+\phi(|Z_s^n|))^2ds+3E\int_0^\tau|Z^n_s|^2ds\\
&\leq&C+ 3TE\int_0^\tau(4\nu^2|Z_s^n|^2+4\nu^2)ds+3E\int_0^\tau|Z^n_s|^2ds\\
&\leq&C+ 3(4\nu^2T+1)E\int_0^\tau|Z_s^n|^2ds.\\
\end{eqnarray*}
Applying It\^{o} formula to $|y^n_t|^2,$ and by (5.4), (5.8), the fact that $\phi$ has a linear growth, and the inequality $2ab\leq \beta a^2+\frac{b^2}{\beta},\ \beta>0,$ we have
\begin{eqnarray*}|y_0^n|^2+E\int_0^\tau|Z^n_s|^2ds&=&E|Y_\tau|^2+2E\int_0^\tau y_s^ng^n_sds
+2E\int_0^\tau y_s^ndA_s^n\\
&\leq& C+2E\int_0^\tau|y_s^n|(\mu|y^n_s|+\phi(|Z_s^n|))ds
+2E\int_0^\tau|y_s^n|dA_s^n\\
&\leq& C+2E\int_0^\tau|y_s^n|(\mu|y^n_s|+\nu|Z_s^n|+\nu)ds
+C[E|A_\tau^n|^2]^\frac{1}{2}\\
&\leq& C+\frac{1}{4}E\int_0^\tau|Z_s^n|^2ds+\frac{1}{6(4\nu^2T+1)}E|A_\tau^n|^2.\\
\end{eqnarray*}
By above two inequalities, we can complete the proof. \ \ $\Box$\\

By (5.4), (5.8), (i) in Proposition 5.6 and linear growth of $\phi$, there exists a constant $C$ independent $n$ such that
$$E\int_0^\tau|g^n_s|^2ds\leq C.\eqno(5.9)$$
By Proposition 5.5, we can get $\forall t\in[0,T],$ there exists $y_{\tau\wedge t}\in L^2({\cal{F}}_{\tau\wedge t}),$ such that
$$y_{\tau\wedge t}^n\rightarrow y_{\tau\wedge t},\ \ \textrm{in}\ \ L^2({\cal{F}}_{\tau\wedge t})\eqno(5.10)$$
as $n\rightarrow \infty$.
By above arguments, we can apply the monotonic limit theorem (see Peng [15, Theorem 2.1] or Peng [16, Theorem 7.2]) to the forward version of  (5.6), then we can get
$$y_{t\wedge\tau}=y_0-A_{t\wedge\tau}-\int^{t\wedge\tau}_0 g_sds+\int^{t\wedge\tau}_0 Z_sdB_s.\ \ t\in[0,T],\eqno(5.11)$$
where $Z_s\in L_{\cal{F}}^2(0,\tau,\textbf{R}^d)$, $g_s\in L_{\cal{F}}^2(0,\tau)$ are the weak limits of $Z_s^n$ and $g_s^n$ in $L_{\cal{F}}^2(0,\tau,\textbf{R}^d)$ and $L_{\cal{F}}^2(0,\tau),$ respectively, $A_t\in{\mathcal{D}}^2_{\cal{F}}(0,\tau)$ is increasing with $A_0=0,$ and for each $t\in[0,T],\ A_{t\wedge\tau}$ is the weak limit of $A_{t\wedge\tau}^n$ in $L^2({\cal{F}}_T).$ By (5.2), Proposition 5.5 and (ii) in Proposition 5.6, we get that as $n\rightarrow\infty,$
$$y^n_{t\wedge\tau}\nearrow Y_{t\wedge\tau},\ \ dP\times dt-a.e.\eqno(5.12)$$
Then by this and Lebesgue dominated convergence theorem, we have
$$y^n\rightarrow Y,\ \ \ \textmd{in}\ \ L_{\cal{F}}^2(0,\tau), \eqno(5.13)$$
Since $y_{t\wedge\tau}$ is RCLL and $Y_{t\wedge\tau}$ is continuous, then by (5.10) and (5.13), we have $\forall t\in[0,T],$
$$y_{t\wedge\tau}=Y_{t\wedge\tau},\ \ P-a.s. \eqno(5.14)$$
Thus $y_{t\wedge\tau}$ is continuous, then by (5.11), we can get $A_t\in {\cal{S}}^2_{{\cal{F}}}(0,\tau)$ and by the monotonic limit theorem in Peng [15, 16] again, we further have
$$Z^n\rightarrow Z,\ \ \ \textmd{in}\ \ L_{\cal{F}}^2(0,\tau), \eqno(5.15)$$
as $n\rightarrow \infty$. By (5.5), (5.13),  (5.15) and the fact that $\phi(|x|)\leq k|x|+\phi(\frac{2\nu}{k})$ for $k\geq2\nu$ (see Fan and Jiang [5, Lemma 4]), we can deduce that the strong limit of $g_t^n$ exists in $L_{\cal{F}}^2(0,\tau)$. Since $g_s\in L_{\cal{F}}^2(0,\tau)$ is the weak limit of $g_s^n$ in $L_{\cal{F}}^2(0,\tau),$ we can get
$$g^n\rightarrow g,\ \ \ \textmd{in}\ \ L_{\cal{F}}^2(0,\tau), \eqno(5.16)$$
as $n\rightarrow \infty$.
Thanks to  (5.10), (5.15) and (5.16), then from (5.6) and (5.11), we can get
$$\forall t\in[0,T],\ A_{\tau\wedge t}^n\rightarrow A_{\tau\wedge t},\ \ \textrm{in}\ \ L^2({\cal{F}}_{\tau\wedge t}),\ \ \ \textrm{and}\ \ \ A^n\rightarrow A,\ \ \textrm{in} \ \ L_{\cal{F}}^2(0,\tau) \eqno(5.17)$$
as $n\rightarrow \infty$. By this and Definition 4.5, we can get that $\forall t\in[0,T],$
$${\cal{E}}_{t\wedge\tau,\tau}[Y_\tau;A^n]\rightarrow{\cal{E}}_{t\wedge\tau,\tau}[Y_\tau;A],\ \ \textrm{in} \ \ L^2({\cal{F}}_T),\eqno(5.18)$$
as $n\rightarrow \infty$. Thus by (5.3), (5.10), (5.14) and (5.18), we have $\forall t\in[0,T],$
$$Y_{t\wedge\tau}={\cal{E}}_{t\wedge\tau,\tau}[Y_\tau;A],\ \ P-a.s.$$
Thanks to (5.10), (5.13)-(5.17), we can complete this proof by passing to limit (a subsequence) of (5.4), (5.6) and (5.7). \ \ $\Box$

%%%%%%%%%%%%%%%%%%%%%%%%%%%%%%%%%%%%%%%%%%%%%%%%%%%%%%%%%%%%%%%%
\section{Representation of ${\cal{F}}$-evaluations by $g$-evaluations}
%%%%%%%%%%%%%%%%%%%%%%%%%%%%%%%%%%%%%%%%%%%%%%%%%%%%%%%%%%%%%%%%
The following representation theorem for ${\cal{F}}$-evaluations is the main result of this paper.\\\\
\textbf{Theorem 6.1} \textit{Let ${\cal{F}}$-evaluation ${\cal{E}}_{s,t}[\cdot]$ satisfy (H1) and (H2).  Then there exists a unique function $g(\omega,t,y,z): \Omega \times [0,T]\times {\mathbf{R}}\times {\mathbf{
R}}^d\longmapsto \mathbf{R},$ satisfying (A1), (A2) and (A3), such that, for each $0\leq s\leq t\leq T$ and $X\in L^2({\cal{F}}_t),$ we have}
\begin{center}
${\cal{E}}_{s,t}[X]={\cal{E}}^g_{s,t}[X],\ \ P-a.s.$
\end{center}
\emph{Proof.}
For $(t,y,z)\in [0,T]\times\textbf{R}\times\textbf{R}^d$, we consider the following process $Y_s^{t,y,z},$ which is the solution of the following SDE on $(t,T]$:

$$dY_s^{t,y,z}=-(\mu|Y_s^{t,y,z}|+\phi(|z|))ds+zdB_s,\ \ \  Y_t^{t,y,z}=y, \eqno(6.1)$$
and the solution of the following BSDE on $[0,t]$:
$$Y_s^{t,y,z}=y+\int_s^t(\mu|Y_r^{t,y,z}|+\phi(|Z_r^{t,y,z}|))dr-\int_s^tZ_r^{t,y,z}dB_r,\ \ s\in[0,t]. \eqno(6.2)$$
Clearly, $Y_s^{t,y,z}\in {\cal{S}}_{\cal{F}}^2(0,T)$ and is an ${\cal{E}}_{s,t}^{\mu,\phi}[\cdot]$-martingale. Then by (i) in Corollary 3.8, we can check that $Y_s^{t,y,z}$ is an ${\cal{E}}_{s,t}[\cdot]$-supermartingale. Now we set the stopping time:
$$\tau_t:=\inf\{s\geq t: |B_s-B_t|\geq1\}\wedge T.\eqno(6.3)$$
Clearly, for $t\in[0,T),$ we have $$|B_{\tau_t}-B_t|=1 \ \textmd{on}\  \{\tau_t<T\},\ \ \textmd{and}\ \ \tau_t>t,\ P-a.s.\eqno(6.4)$$
By (6.1) and (6.3), we have for $s\in[t,T],$
$$|Y_{s\wedge\tau_t}^{t,y,z}|\leq |y|+\int_t^{s\wedge\tau_t}\mu|Y_r^{t,y,z}|dr+\phi(|z|)T+|z|,\ \ P-a.s.$$
Then by Gronwall's inequality, we can get for $s\in[t,T],$
$$|Y_{s\wedge\tau_t}^{t,y,z}|\leq( |y|+|z|+\phi(|z|)T)e^{\mu T},\  \  P-a.s. \eqno(6.5)$$
By (6.2), Lemma 2.6 and (6.5), we have $Y_{s}^{t,y,z}\in {\cal{S}}_{\cal{F}}^\infty(0,\tau_t).$
Then by Theorem 5.4, there exists a process $A_s^{t,y,z}\in {\cal{S}}^2_{{\cal{F}}}(0,\tau_t)$, which is increasing with $A_0^{t,y,z}=0$ such that $\forall s\in[0,T],$
$${\cal{E}}_{s\wedge{\tau_t},{\tau_t}}[Y_{{\tau_t}}^{t,y,z};A^{t,y,z}]
=Y_{s\wedge{\tau_t}}^{t,y,z},\ \ P-a.s.,$$
and there exists a pair $(g_r^{t,y,z}, Z^{t,y,z}_r)$ such that
$$Y_{s\wedge{\tau_t}}^{t,y,z}=Y_{{\tau_t}}^{t,y,z}+A_{{\tau_t}}^{t,y,z}
-A_{s\wedge{\tau_t}}^{t,y,z}+\int_{s\wedge{\tau_t}}^{{\tau_t}}g_r^{t,y,z}dr
-\int_{s\wedge{\tau_t}}^{{\tau_t}}Z^{t,y,z}_rdB_r,\ \ P-a.s.,\ \ s\in[0,T],\eqno(6.6)$$
$$|g_{s}^{t,y,z}|\leq\mu|Y_{s}^{t,y,z}|+\phi(|Z^{t,y,z}_{s}|),\ \ dP\times dt-a.e.,\ \ s\in[0,\tau_t],\eqno(6.7)$$
and for $(t',y',z')\in [0,T]\times\textbf{R}\times\textbf{R}^d$,
$$|g_{s}^{t,y,z}
-g_{s}^{t',y',z'}|\leq \mu|Y_{s}^{t,y,z}-Y_{s}^{t',y',z'}|
+\phi(|Z_{s}^{t,y,z}
-Z_{s}^{t',y',z'}|),\ \ dP\times dt-a.e.,\ \  s\in[0,{\tau_t}\wedge{\tau_{t'}}]. \eqno(6.8)$$
For each $t''\geq t$ and $X\in L^\infty({\cal{F}}_{t''}),$  we set
$$\bar{Y}_s^{t'',X}:={\cal{E}}_{s,t''}[X].$$
By Theorem 5.4, there exists a pair $(\bar{g}_r^{t'',X}, \bar{Z}^{t'',X}_r)$ such that
$$\bar{Y}_s^{t'',X}=X+\int_s^{t''}{\bar{g}}_r^{t'',X}dr-\int_s^{t''}\bar{Z}^{t'',X}_rdB_r,\ \ s\in[0,t''].\eqno(6.9)$$
and
$$|{g}_{s}^{t,y,z}-\bar{g}_{s}^{t'',X}|\leq \mu|{Y}_{s}^{t,y,z}-\bar{Y}_{s}^{t'',X}|
+\phi(|{Z}^{t,y,z}_{s}-\bar{Z}^{t'',X}_{s}|),\ \ dP\times dt-a.e.,\ s\in[0,{\tau_t}\wedge t''].\eqno(6.10)$$
Comparing the bounded variation parts and martingale parts and of (6.1) and (6.6), we get
$$Z^{t,y,z}_{s}= z,\ \ s\in[t,\tau_t],\ \ dP\times dt-a.e.$$
From this,  we can rewrite (6.7), (6.8) and (6.10) as
$$|g_{s}^{t,y,z}|\leq\mu|Y_{s}^{t,y,z}|+\phi(|z|),\ \ dP\times dt-a.e.,\ \ s\in[t,\tau_t],\eqno(6.11)$$
$$|g_{s}^{t,y,z}-g_{s}^{t',y',z'}|
\leq\mu|Y_{s}^{t,y,z}-Y_{s}^{t',y',z'}|
+\phi(|z-z'|),\ \ dP\times dt-a.e.,\ \ s\in[t\vee t',{\tau_t}\wedge{\tau_{t'}}],\eqno(6.12)$$
and
$$|{g}_{s}^{t,y,z}-\bar{g}_{s}^{t'',X}|\leq \mu|{Y}_{s}^{t,y,z}-\bar{Y}_{s}^{t'',X}|+\phi(|z-\bar{Z}^{t'',X}_{s}|),\ \ dP\times dt-a.e.,\ \ s\in[t,{\tau_t}\wedge t''],\eqno(6.13)$$
respectively. Now for $n\geq1,$ we set $t^n_i=i2^{-n}T,\ i=0,1,2\cdots,2^n,$ and
$$g^n(s,y,z):=\sum_{i=0}^{2^n-1}g_{s}^{t^n_i,y,z}1_{[t^n_i,\tau_{t^n_i}\wedge t^n_{i+1})}(s),\ \ \textrm{for}\  (s,y,z)\in [0,T)\times\textbf{R}\times\textbf{R}^d.$$
Clearly, for each $n\geq1$ and each $s\in[0,T),$ there always exists an interval denoted by $[t_{i_s}^n,t_{{i_s}+1}^n),$ such that $s\in[t_{i_s}^n,t_{{i_s}+1}^n).$ Thus we have
$$g^n(s,y,z)=g_{s}^{t^n_{i_s},y,z}1_{\{s<\tau_{t^n_{i_s}}\}},\ \ \textrm{for}\  (s,y,z)\in [0,T)\times\textbf{R}\times\textbf{R}^d.\eqno(6.14)$$
By (6.14), (6.11) and (6.5), there exists a constant $C$ only dependent on $y, z, \mu, \nu$ and $T$ such that
$$\|g^n(s,y,z)\|_{L_{\cal{F}}^\infty(0,T)}\leq C.\eqno(6.15)$$ Moreover, we have\\\\
\
\textbf{Proposition 6.2}\textit{ For $(s,y,z)\in [0,T]\times\textbf{R}\times\textbf{R}^d,$ $g^n(s,y,z)$ is a Cauchy sequence in $ L^2_{{\cal{F}}}(0,T).$}\\\\
\emph{Proof.} For $(s,y,z)\in [0,T)\times\textbf{R}\times\textbf{R}^d,$ by (6.1) and the classic estimate on solutions of SDEs, we have
\begin{eqnarray*}
\ \ \    \ \ \ \    \ \ \ \ \    E\left[|Y_{s}^{t^n_{i_s},y,z}-y|^2\right]
&\leq&E\left|\int_{t^n_{i_s}}^{s}(\mu|Y_{r}^{t^n_{i_s},y,z}|
+\phi(|z|))dr+z(B_{s}-B_{t^n_{i_s}})\right|^2\\
&\leq& 2^{-n}C(|y|^2+|z|^2+1),\ \ \ \ \ \ \ \ \ \ \  \ \ \ \ \ \  \ \ \ \ \ \ \ \  \ \ \ \ \ \ \   \ \      \ \ \ \    \ \ \ \ \ \ \  (6.16)
\end{eqnarray*}
where $C$ is a constant only dependent on $\mu$, $\nu$ and $T$.

For $s\in[0,T),$ we set $\underline{\tau}_s:=\liminf_{n\rightarrow\infty}\tau_{t^n_{i_s}}.$ Clearly, $\underline{\tau}_s$ is a stopping time, and we can get for a.e. $\omega\in\Omega,$ there exists a sequence $\{\tau_{t_{i_s}^{n_\omega}}\}_{n_\omega\geq1}$ such that $\underline{\tau}_s(\omega)=\lim_{n_\omega\rightarrow\infty}\tau_{t^{n_\omega}_{i_s}}(\omega).$  By this and (6.4), we can further have  for a.e. $\omega\in\Omega,$
$$|B_{\underline{\tau}_s(\omega)}(\omega)-B_s(\omega)|
=\lim_{n\rightarrow\infty}|B_{\tau_{t^{n_\omega}_{i_s}}(\omega)}(\omega)-B_{t_{i_s}^{n_\omega}} (\omega)|=1,\ \ \textmd{if}\ \ \underline{\tau}_s(\omega)<T.$$
From this, (6.3) and (6.4), it follows that for each $s\in[0,T),$
$$\underline{\tau}_s\geq\tau_s>s,\ \ P-a.s.$$
Thus, for two integers $m, n$ and any $\varepsilon>0,$  we have for each $s\in[0,T),$
\begin{eqnarray*}
&&\lim_{m,n\rightarrow\infty}P\left(1_{\{s\geq{\tau_{t_{i_s}^n}}\wedge\tau_{t_{i_s}^m}\}}|g^n(s,y,z)-g^m(s,y,z)|^2>\varepsilon\right)\\
&\leq&\lim_{m,n\rightarrow\infty}P \left(s\geq{\tau_{t_{i_s}^n}}\wedge\tau_{t_{i_s}^m}\right)\\
&\leq&\lim_{m,n\rightarrow\infty}P \left(s\geq\inf_{k\geq n}\tau_{t_{i_s}^k}\wedge\inf_{l\geq m}\tau_{t_{i_s}^l}\right)\\
&=&P \left(\cap_{m,n\geq1}\left\{s\geq\inf_{k\geq n}\tau_{t_{i_s}^k}\wedge\inf_{l\geq m}\tau_{t_{i_s}^l}\right\}\right)\\
&=&P \left(s\geq\underline{\tau}_s\right)\\
&=&0.
\end{eqnarray*}
By this, (6.15) and dominated convergence theorem, we have for each $s\in[0,T),$
$$\lim_{m,n\rightarrow\infty}E\left[1_{\{s\geq{\tau_{t_{i_s}^n}}\wedge\tau_{t_{i_s}^m}\}}
|g^n(s,y,z)-g^m(s,y,z)|^2\right]=0.\eqno(6.17)$$
By (6.14), (6.12) and (6.16), we have for $a.e., s\in[0,T],$
\begin{eqnarray*}
\ \ \ \  \ \ &&E\left[1_{\{s<{\tau_{t_{i_s}^n}}\wedge\tau_{t_{i_s}^m}\}}|g^n(s,y,z)-g^m(s,y,z)|^2\right]\\
&=& E\left[1_{\{s<{\tau_{t_{i_s}^n}}\wedge\tau_{t_{i_s}^m}\}}|g_{s}^{t^n_{i_s},y,z}
-g_{s}^{t^m_{i_s},y,z}|^2\right]\\
&=& E\left[1_{\{t_{i_s}^n\vee t_{i_s}^m\leq s<{\tau_{t_{i_s}^n}}\wedge\tau_{t_{i_s}^m}\}}|g_{s}^{t^n_{i_s},y,z}
-g_{s}^{t^m_{i_s},y,z}|^2\right]\\
&\leq& E\left[\mu^2|Y_{s}^{t^n_{i_s},y,z}
-Y_{s}^{t^m_{i_s},y,z}|^2\right]\\
&\leq& 2E\left[\mu^2(|Y_{s}^{t^n_{i_s},y,z}
-y|^2+|Y_{s}^{t^m_{i_s},y,z}-y|^2)\right]\\
&\leq&2\mu^2\left(2^{-n}C(|y|^2+|z|^2+1)+2^{-m}C(|y|^2+|z|^2+1)\right). \ \ \  \ \ \ \  \ \ \  \ \ \ \ \      \ \ \ \ \ \ \ \ \ \ \ \ \ \ \ \ \ \  (6.18)
\end{eqnarray*}
By (6.17) and (6.18), we have for $a.e., s\in[0,T],$
\begin{eqnarray*}
&&\lim_{m,n\rightarrow\infty}E\left[|g^n(s,y,z)-g^m(s,y,z)|^2\right]\\
&\leq&\lim_{m,n\rightarrow\infty}E\left[1_{\{s<{\tau_{t_{i_s}^n}}\wedge\tau_{t_{i_s}^m}\}}|g^n(s,y,z)-g^m(s,y,z)|^2\right]\\
&&+\lim_{m,n\rightarrow\infty}E\left[1_{\{s\geq{\tau_{t_{i_s}^n}}\wedge\tau_{t_{i_s}^m}\}}|g^n(s,y,z)-g^m(s,y,z)|^2\right]\\
&=&0.
\end{eqnarray*}
By this, Fubini's Theorem, (6.15) and dominated convergence theorem, we have
\begin{eqnarray*}
&&\lim_{m,n\rightarrow\infty}E\int_0^T|g^n(s,y,z)-g^m(s,y,z)|^2ds
\\&\leq&\lim_{m,n\rightarrow\infty}\int_0^T E|g^n(s,y,z)-g^m(s,y,z)|^2ds\\
&=&0.
\end{eqnarray*}
The proof is complete.\ \ $\Box$\\

We denote the limit of $g^n(s,y,z)$ in $ L^2_{{\cal{F}}}(0,T)$ by $g(s,y,z).$ We can further get the following properties.\\\\
\
\textbf{Proposition 6.3}  \textit{ $g(s,y,z)$ satisfies (A1)-(A3) and for $a.e., s\in[0,t''],$
$$|{g}(s,y,z)-{\bar{g}}_s^{t'',X}|\leq \mu|y-\bar{Y}_s^{t'',X}|+\phi(|z-\bar{Z}^{t'',X}_s|),\ \ P-a.s.\eqno(6.19)$$}
\emph{Proof. } By (6.15), we have $g(s,y,z)$ satisfies (A2). By (6.14), (6.11) and (6.5), we have $g^n(t,0,0)=0,$ $dP\times dt-a.e.$ Thus  $g(s,y,z)$ satisfies (A3). By (6.14) and (6.12), we can get $dP\times dt-a.e.,$
\begin{eqnarray*}
&&|g^n(s,y,z)-g^n(s,y',z')|\\
&=&1_{ \{ s<\tau_{t_{i_s}^n}\}}|g_{s}^{t^n_{i_s},y,z}
-g_{s}^{t^n_{i_s},y',z'}|\\
&\leq&\mu|Y_{s}^{t^n_{i_s},y,z}
-Y_{s}^{t^n_{i_s},y',z'}|+\phi(|z-z'|)\\
&\leq&\mu\left(|Y_{s}^{t^n_{i_s},y,z}-y|
+|Y_{s}^{t^n_{i_s},y',z'}-y'|\right)+\mu|y-y'|+\phi(|z-z'|).
\end{eqnarray*}
Then from Proposition 6.2 and (6.16), it follows that $g(s,y,z)$ satisfies (A1). By (6.14) and (6.13), we have for $a.e., s\in[0,t''],$ $P-a.s.,$
\begin{eqnarray*}
&&|g^n(s,y,z)-{\bar{g}}_s^{t'',X}|\\
&=& 1_{ \{ s<\tau_{t_{i_s}^n}\}}|g^n(s,y,z)-{\bar{g}}_s^{t'',X}|+1_{ \{ s\geq\tau_{t_{i_s}^n}\}}|g^n(s,y,z)-{\bar{g}}_s^{t'',X}|\\
&=&  1_{ \{ s<\tau_{t_{i_s}^n}\}}|g_{s}^{t^n_{i_s},y,z}-{\bar{g}}_{s}^{t'',X}|+1_{ \{ s\geq\tau_{t_{i_s}^n}\}}|g^n(s,y,z)-{\bar{g}}_s^{t'',X}|
\\&\leq&\left(\mu|Y_{s}^{t^n_{i_s},y,z}-\bar{Y}_{s}^{t'',X}|
+\phi(|z-\bar{Z}^{t'',X}_{s}|)\right)+1_{ \{ s\geq\tau_{t_{i_s}^n}\}}|g^n(s,y,z)-{\bar{g}}_s^{t'',X}|\\
\\&\leq&\left(\mu |Y_{s}^{t^n_{i_s},y,z}-y|
+\mu|y-\bar{Y}_{s}^{t'',X}|+\phi(|z-\bar{Z}^{t'',X}_{s}|)\right)+1_{ \{ s\geq\tau_{t_{i_s}^n}\}}|g^n(s,y,z)-{\bar{g}}_s^{t'',X}|, \ \
\end{eqnarray*}
By Proposition 6.2, (6.16) and the argument of (6.17), we can obtain (6.19). $\Box$\\

Now, we come back the proof of Theorem 6.1. For fixed $t\in[0,T]$ and $X\in L^\infty({\cal{F}}_t),$  we set
$$\bar{Y}_s^{t,X}:={\cal{E}}_{s,t}[X],\ s\in[0,t].$$
Then by Theorem 5.4,  there exists a pair $(\bar{g}_u^{t,X}, \bar{Z}^{t,X}_u)$ such that for $s\in[0,t],$
$$\bar{Y}_s^{t,X}=X+\int_s^{t}\bar{g}_u^{t,X}du-\int_s^{t}\bar{Z}^{t,X}_udB_u.$$
We consider the following BSDE on $[0,t],$
$${Y}_s^{t,X}=X+\int_s^{t}g(u,{Y}_u^{t,X},{Z}_u^{t,X})du-\int_s^{t}{Z}^{u,X}_udB_u.$$
Set $\hat{g}_s:=g(s,{Y}_s^{t,X},{Z}_s^{t,X})-\bar{g}_s^{t,X},$ ${\hat{Y}}_s:={Y}_s^{t,X}-\bar{Y}_s^{t,X}$ and ${\hat{Z}}_s:={Z}_s^{t,X}-\bar{Z}_s^{t,X}.$ By (6.19) and (2.10), we have for $s\in[0,t]$
$$|\hat{g}_s|\leq \mu|{\hat{Y}}_s|+\phi(|{\hat{Z}}_s|)\leq \mu|{\hat{Y}}_s|+n|{\hat{Z}}_s|+\phi\left(\frac{2\nu}{n}\right),\ \ dP\times dt-a.e.,\ \ \textrm{ for}\ \ n\geq2\nu.$$
By this and the proof of uniqueness of solutions of BSDEs in Fan and Jiang [5, Theorem 2], we can get $\forall s\in[0,t],$ $P-a.s.,$ ${Y}_s^{t,X}=\bar{Y}_s^{t,X}.$ For $X\in L^2({\cal{F}}_t),$ we set $X_n=(X\vee(-n))\wedge n.$ Thus, we have ${\cal{E}}_{s,t}[X_n]={\cal{E}}^g_{s,t}[X_n].$ By this, Lemma 2.5 and Lemma 3.9, we have $\forall s\in[0,t],$
$${\cal{E}}_{s,t}[X]={\cal{E}}^g_{s,t}[X],\ P-a.s.$$

Now, we prove the uniqueness of $g$. Suppose there exists  another function $\bar{g}(\omega,t,y,z): \Omega \times [0,T]\times {\textbf{R}} \times{\mathbf{
R}}^d\longmapsto \mathbf{R}$ satisfying (A1), (A2) and (A3), such that for each $t\in[0,T],$ $X\in L^2({\cal{F}}_t),$ we have for all $s\in[0,t],$ ${\cal{E}}^g_{s,t}[X]={\cal{E}}^{\bar{g}}_{s,t}[X],\ P-a.s.$ Then as the argument in the proof of Zheng and Li [19, Theorem 5.1], we can get $dP\times dt-a.e.,$
$$g(t,y,z)=\bar{g}(t,y,z),\ \forall(y,z)\in {\textbf{R}} \times{\mathbf{
R}}^d,$$
from the representation theorem for generators of BSDEs (see Fan and Jiang [4, Theorem 2] or Jia [9, Theorem 3.4]). The proof is complete. \ \ $\Box$\\\\
\
\textbf{Corollary 6.4}\textit{ Let ${\cal{F}}$-evaluation ${\cal{E}}_{s,t}[\cdot]$ satisfy (H1) and (H2), $K\in {\cal{D}}^2_{{\cal{F}}}(0,T).$  Then there exists a unique function $g(\omega,t,y,z): \Omega \times [0,T]\times {\mathbf{R}}\times {\mathbf{
R}}^d\longmapsto \mathbf{R},$ satisfying (A1), (A2) and (A3), such that, for each $0\leq s\leq t\leq T$ and $X\in L^2({\cal{F}}_t),$ we have}
$${\cal{E}}_{s,t}[X;K]={\cal{E}}^g_{s,t}[X;K],\ \ P-a.s.\eqno(6.20)$$
\emph{Proof.} We sketch this proof. By Theorem 6.1 and Proposition 3.5, we can get there exists a unique function $g(\omega,t,y,z): \Omega \times [0,T]\times {\mathbf{R}}\times {\mathbf{
R}}^d\longmapsto \mathbf{R},$ satisfying (A1), (A2) and (A3), such that, for each $K\in {\cal{D}}^{2,0}_{{\cal{F}}}(0,T),$ we have (6.20).
Thus, for $K\in {\cal{D}}^2_{{\cal{F}}}(0,T),$ by Definition of ${\cal{E}}_{s,t}[X;K]$ and Lemma 2.5, we can still get (6.20). The proof is complete.\ \ $\Box$\\\\
\textbf{Remark 6.5}
\begin{itemize}
\item[(i)] Theorem 5.1 and Theorem 5.2 are existence and uniqueness theorem and comparison theorem
of ${\cal{E}}(f,X,T)$, respectively, with $X\in L^\infty({\cal{F}}_T)$ and $f(\cdot,0)\in L^\infty_{{\cal{F}}}(0,T).$  By Corollary 6.4 and the similarly argument as Zheng and Li [19, Corollary 5.1], we can get that the two theorems are both true for ${\cal{E}}(f,X,T)$ with $X\in L^2({\cal{F}}_T)$ and $f(\cdot,0)\in L^2_{{\cal{F}}}(0,T).$
\item[(ii)] In Theorem 6.1, if ${\cal{E}}_{s,t}^{\mu,\phi}[\cdot]$ is placed by ${\cal{E}}_{s,t}^{\mu,\mu}[\cdot],$ then Theorem 6.1 will become Peng [14, Theorem 3.1]. In Theorem 6.1, if the ${\cal{F}}$-evaluation become an ${\cal{F}}$-expectation, then (H1) will become (H1) in Zheng and Li [19], and by Zheng and Li [19, Remark 3.1], the ${\cal{F}}$-evaluation will satisfy translation invariance ((H2) in Zheng and Li [19]). By this, we can further get that $g$ in Theorem 6.1 will be independent on $y$ (see Jia [8, Corollary 2.3.14]). Thus Theorem 6.1 will become Zheng and Li [19, Theorem 5.1].
\item[(iii)] In Theorem 6.1, can we replace the domination condition (H1) by the following (H4)?

(H4) : For each $0\leq s\leq t\leq T$ and $X,\ Y$ in $L^2({\cal{F}}_t),$ we have
\begin{center}
${\cal{E}}_{s,t}[X]-{\cal{E}}_{s,t}[Y]\leq{\cal{E}}_{s,t}^{\phi_1,\phi_2}[X-Y],\ \ P-a.s.$
\end{center}
where $\phi_1(\cdot)$ and $\phi_2(\cdot)$ are functions given in (A1).

In general, the solution of BSDE with generator $g=\phi_1(|y|)+\phi_2(|z|),$ denoted by ${\cal{E}}_{s,t}^{\phi_1,\phi_2}[\cdot],$ is not unique (see Jia [8, Remark 3.2.5]). Consequently, under (H4), we can not obtain a representation theorem like Theorem 6.1 using the method in this paper.
\end{itemize}
\
\\
\textbf{Acknowledgements}\ \   The authors are supported by the National Natural Science Foundation of China (No. 11571024) and Natural Science Foundation of Beijing (No. 1132008). The authors would like to thank the anonymous referee and Associate Editor for their valuable comments and suggestions
which improved the presentation of this paper. The authors would also like to thank Professor Renming Song for his help and useful suggestions. The first author would like to thank Professor Zengjing Chen for providing him a good opportunity for studying at Qilu Institute of Finance of Shandong University Spring 2015, where part of this paper was written.

\end{document}